\documentclass[times, 11pt, dvips]{article}
\usepackage{amssymb, amsthm, hyperref, xypic}
\hypersetup{pdftitle={Gluing Seiberg-Witten Monopoles}, verbose,
            hyperindex=true,
            pdfpagemode=UseThumbs,       %% options: None, FullScreen
            colorlinks=true,
            linkcolor=blue,
            citecolor=magenta,
            urlcolor=red,
            pagecolor=green,
            menucolor=cyan,
            anchorcolor=black,
            naturalnames=true,
            pdfstartview=XYZ
}
\makeglossary
\makeindex
\input{boxedeps.tex}
\SetepsfEPSFSpecial
\HideDisplacementBoxes
\begin{document}
%
%%%%%%%%%%%%%%%% Layout %%%%%%%%%%%%%%%%%
%
% \pagenumbering{arabic}                 %% or roman
% \parskip=2cm
% \parindent=0cm
% \baselineskip=24pt
% \renewcommand{\baselinestretch}{1.6}
%
%%%%%%%%%%%%%% Definitions %%%%%%%%%%%%%%
%
% \input{def}
%
\def\sw{Seiberg-Witten }
\def\SW{\ifmmode{\mathrm{SW}}\else{Seiberg-Witten }\fi}
\def\mfld{manifold }
\def\mflds{manifolds }
\def\4mfld{four-manifold}
\def\3mfld{three-manifold}
\def\nbhd{neighborhood }
\def\wrt{with respect to }
\def\cend{cylindrical end }
\def\cyend{cylindrical-end }
\def\cancel#1#2{\ooalign{$\hfil#1\mkern1mu/\hfil$\crcr$#1#2$}}
\def\dirac{\mathpalette\cancel\partial}
\def\cinf{C^{\infty}}
\def\spin{\ifmmode{\mathrm{Spin}}\else{Spin}\fi}
\def\spinc{\ifmmode{\mathrm{Spin}^c}\else{Spin$^c$}\fi}
\def\Cl{\ifmmode{\bbC{\mathrm{l}}}\else{$\bbC{\mbox{l}}$}\fi}
\def\sminus{\setminus}
\def\plus{\oplus}
\def\diff{\simeq}
\def\iso{\cong}
\def\hom{\mathrm{Hom}}
\def\mend{\mathrm{End}}
\def\aut{\mathrm{Aut}}
\def\im{\mathrm{Im}}
\def\ker{\mathrm{Ker}}
\def\bbC{\mbox{${\mathbb C}$}}
\def\bbR{\mbox{${\mathbb R}$}}
\def\bbQ{\mbox{${\mathbb Q}$}}
\def\bbZ{\mbox{${\mathbb Z}$}}
\def\A{{\cal A}}
\def\B{{\cal B}}
\def\C{{\cal C}}
\def\D{{\cal D}}
\def\E{{\cal E}}
\def\F{{\cal F}}
\def\G{\ifmmode{\cal G}\else{${\cal G}$}\fi}
\def\M{\ifmmode{\cal M}\else{${\cal M}$}\fi}
\def\N{\ifmmode{\cal N}\else{${\cal N}$}\fi}
\def\H{{\cal H}}
\def\R{{\cal R}}
\def\L{{\cal L}}
\def\W{{\cal W}}
\def\S{{\cal S}}
\def\c{{\bf c}}
\def\cfrak{\ifmmode{\mathfrak c}\else{${\mathfrak c }$}\fi}
\def\s{\ifmmode{\mathfrak s}\else{${\mathfrak s }$}\fi}
\def\t{\ifmmode{\mathfrak t}\else{${\mathfrak t }$}\fi}
\def\m{\ifmmode{\frak m}\else{${\frak m }$}\fi}
\def\gd{{\frak D}}
\def\xibold{\ifmmode{\mbox{\boldmath$\xi$\unboldmath}}\else{\boldmath$\xi$\unboldmath}\fi}
\def\xipt{\xibold}
\def\xivector{\xi}
\def\xihat{\hat\xivector}
\def\etabold{\ifmmode{\mbox{\boldmath$\eta$\unboldmath}}\else{\boldmath$\eta$\unboldmath}\fi}
\def\etapt{\etabold}
\def\etavector{\eta}
\def\zetabold{\ifmmode{\mbox{\boldmath$\zeta$\unboldmath}}\else{\boldmath$\zeta$\unboldmath}\fi}
\def\zetapt{\zetabold}
\def\zetavector{\zeta}
\def\zetahat{\hat\zetavector}
\def\wedgedot{\,\dot{\wedge}\,}
\def\ip#1#2{\langle #1,#2\rangle}        %% inner product
\def\iint{\mathop{\int\!\int}}           %% double integral
\def\inte#1#2{\int\limits_{#1}^{#2}}
\def\sl#1#2#3{\sum\limits_{#1=#2}^{#3}}
\def\forms#1#2{\Omega^#1_#2}             %% #1-forms in L2_#2 norm
\def\wforms#1#2#3{\Omega^#1_{#2,#3}}     %% weighted forms ...
%
%% My defined colors, using the graphics package:
%
\definecolor{brown}{rgb}{0.38,0.19,0}
\definecolor{mygreen}{rgb}{0.13,0.38,0.13}
%
%% New CD macros (can just be \def'ed. (?))
%
\newcommand{\cdr}{\longrightarrow}
\newcommand{\cdl}{\longleftarrow}
\newcommand{\cdu}{\uparrow}
\newcommand{\cdd}{\downarrow}
\newcommand{\cdur}{\nearrow}   %% up right arrow
\newcommand{\cdul}{\nwarrow}   %% up left arrow
\newcommand{\cddr}{\searrow}   %% down rigth arrow
\newcommand{\cddl}{\swarrow}   %% down left arrow
%\newcommand{\cdrl}{\cd\rightleftarrows}
%\newcommand{\cdlr}{\cd\leftrightarrows}
%\newcommand{\cdud}{\cd\updownarrow}
%\newcommand{\cddu}{\cd\downuparrow}
%
%% Environment definitions
%
%% \newtheorem{type}{title}[in_counter]  % counter resets as the
                                         % in_counter increases. 
%% \newtheorem{type}[num_like]{title}    % will be numbered in the
                                         % num_like series.
\newtheorem{thm}{Theorem}
\newtheorem{theorem}{Theorem}[section]
\newtheorem{prop}[theorem]{Proposition}
\newtheorem{lemma}[theorem]{Lemma}
\newtheorem{claim}[theorem]{Claim}
\newtheorem{ex}[theorem]{Example}
\newtheorem{cor}[theorem]{Corollary}
\newtheorem{Def}[theorem]{Definition}
\newtheorem{addendum}[theorem]{Addendum}
\newtheorem{rem}[theorem]{Remark}
%
%% \theoremstyle{style}        % amsthm command -- can be issued at any
                               % point in the text and affects the sequel. 
                               % style options: plain, definition, remark
\theoremstyle{definition}                % amsthm command
\newtheorem{statement}{Statement}[section]
\newtheorem{st}[statement]{}
% \numberwithin{equation}{section}         % amsmath command
\theoremstyle{plain}                     % amsthm command
\newtheorem*{thm*}{Theorem}              % amsthm command
\def\pf{\vspace{2ex}\noindent{\bf Proof.} }
\def\tpf#1{\vspace{2ex}\noindent{\bf Proof of Theorem #1.} }
\def\ppf#1{\vspace{2ex}\noindent{\bf Proof of Proposition #1.} }
\def\lpf#1{\vspace{2ex}\noindent{\bf Proof of Lemma #1.} }
\def\cpf#1{\vspace{2ex}\noindent{\bf Proof of Corollary #1.} }
\def\clpf#1{\vspace{2ex}\noindent{\bf Proof of Claim #1.} }
\newenvironment{Defn}{\begin{Def} \rm }{\end{Def}}
\newenvironment{remark}{\begin{rem}\rm }{\end{rem}}
\newenvironment{squote}{\begin{quote}\small}{\end{quote}}
\newenvironment{Proof}{\begin{squote}\pf}{\endpf\end{squote}}
\newenvironment{tproof}[1]{\begin{squote}\tpf#1}{\endpf\end{squote}}
\newenvironment{pproof}[1]{\begin{squote}\ppf#1}{\endpf\end{squote}}
\newenvironment{lproof}[1]{\begin{squote}\lpf#1}{\endpf\end{squote}}
\newenvironment{cproof}[1]{\begin{squote}\cpf#1}{\endpf\end{squote}}
\newenvironment{clproof}[1]{\begin{squote}\clpf#1}{\endpf\end{squote}}
\newenvironment{comments}{\smallskip\noindent{\bf Comments:}\begin{enumerate}}{\end{enumerate}\smallskip}
\newenvironment{matris}[4]{\left(\begin{array}{cc}#1&#2\\#3&#4}{\end{array}\right)}
\newenvironment{bordar}[2]{\left(\begin{array}{c}#1\\#2}{\end{array}\right)}
\newcommand{\glossaryentry}[3]{\parbox[t]{3cm}{#1}\parbox[t]{10.5cm}{#2}\hfill#3\\}
\newcommand{\note}[1]{\marginpar{\scriptsize #1 }}
%
% End of Proof Symbol at the end of an equation must precede $$.
% (I need to insert another pair of dollars here to fix their parity
% in the sequel of the text.) $$
\def\endpf{\relax\ifmmode\expandafter\endproofmath\else
  \unskip\nobreak\hfil\penalty50\hskip.75em\hbox{}\nobreak\hfil\bull
  {\parfillskip=0pt \finalhyphendemerits=0 \bigbreak}\fi}
\def\endproofmath$${\eqno\bull$$\bigbreak}
\def\bull{\vbox{\hrule\hbox{\vrule\kern3pt\vbox{\kern6pt}\kern3pt\vrule}\hrule}}
%
%%%%%%%%%%%%%%%%%%%%%%%%%%% Begin Text Here %%%%%%%%%%%%%%%%%%%%%%%%%
%
\title{Gluing Seiberg-Witten Monopoles
\thanks{\href{http://www.ams.org/msc/}{AMS Mathematics Subject Classification} 2000: Primary 57R57; Secondary 58D27, 57R58.}
}
\author{\hyperlink{myaddress}{Pedram Safari}
% \thanks{The author is grateful to John Morgan.}
}
\date{}     % or: \date{\today}
% November 19, 2003
\maketitle
\thispagestyle{empty}
\begin{abstract}
  We establish a canonical gluing procedure for Seiberg-Witten
  monopoles on the two pieces of a closed, oriented \4mfld $X$ which
  is split along a 3-dimensional closed, oriented submanifold. We only
  assume that the (unperturbed) character variety is Kuranishi-smooth
  and the limiting maps are transversal --- then we will be able to
  glue regular monopoles over the irreducible points of the character
  variety.  
% This paper deals with the infinitesimal aspect of the gluing.
\end{abstract}
\tableofcontents
\newpage
\section{Introduction}
\label{sec:intro}%
% We are going to give a brief account of the history of the theory and
% an outline of the proof.
% \par
The advent of \sw theory in 1994 led not only to a great
simplification of the gauge-theoretic results obtained earlier by
Donaldson, but also to new advances such as proofs for the Thom
conjecture \cite{km}. One of the advantages of the \sw theory is that
the bubbling phenomenon does not occur, thus resulting in a compact
moduli space for closed \4mflds.
% Also, the gauge group is abelian, thus eliminating much of the
% technical difficulties.
\par
Naturally, finding methods to compute these invariants would be
desirable. For symplectic manifolds, Taubes settled this question by
relating the invariants to those of Gromov-Witten \cite{tswgr},
\cite{tswgr1}. Fintushel and Stern
% \note{is this precise?}%
described how the \sw invariants change under certain surgeries over a
knot \cite{fs}. In yet another direction, one could ask if it is
possible to compute the \sw invariants of a \4mfld which is decomposed
into two parts, given the relevant information on the pieces.
\par
To make this more precise, let $X$ be a closed, oriented \4mfld, $Y$ a
closed, oriented, embedded, dividing submanifold, and $X^+$ and $X^-$
the two components of $X\setminus Y.$ $X^\pm$ could be considered as
\cyend manifolds with ends isometric to $I\times Y,$ where $I$ is an
interval. If $X^\pm$ are equipped with \spinc-structures agreeing on
$Y,$ then these structures give rise to a unique \spinc-structure on
$X.$ We are interested in the \sw moduli space of $X$ in terms of
those of $X^+,$ $X^-$ and $Y.$
\par
There has already been several partial results around this question; see
\cite{mst}, \cite{moy}, \cite{joe}, \cite{mcw1}, \cite{mw2}, \cite{c}. The
main concern in those works is to explicitly find the solutions to the
Seiberg-Witten equations by analytical means, mostly dealing with some
particular type of the boundary manifold $Y,$ such as circle bundles
over Riemann surfaces or Seifert fibered spaces.
% \index{circle bundles}
% \index{Riemann surfaces}
\index{Seifert fibered space}%
\par
Here we % take on a different road to 
develop a more general gluing scheme for \sw monopoles, thus setting
the ground for a Mayer-Vietoris type theorem for Seiberg-Witten moduli
spaces; this would be the subject of a forthcoming paper. Our
approach is essentially based on Taubes' method for constructing
glued-up ASD connections on connected sums, and is an adaptation of
the work of Morgan and Mrowka in gluing ASD SU(2)-connections
(\cite{mm}) in the context of \sw theory.
% \par
% \ldots
\par
We will assume that the \3mfld $Y$ has no reducible point in its \sw
moduli space, otherwise we would simply confine ourselves to the open
set of irreducible points. % in $\M(Y)$. 
This could also be achieved by a small perturbation in the \sw
equations, which would relieve us from singular points as well,
producing a zero-dimensional moduli space, but we are not heeding much
this way. Instead, we would allow a positive-dimensional, ``smooth''
(in the sense of Kuranishi) moduli space, which is free from
reducibles, while the obstruction spaces can be non-trivial. On
$X^\pm,$ we would require the stronger condition of regularity for
monopoles.
\par
There are practical motives for this unperturbative approach --- in
many concrete examples, such as the case of Seifert fibered spaces, we
explicitly know the solutions to the {\it unperturbed} equations on
the three-manifold. However, our method should well work in the
perturbative setup as well.
% \note{?}%
% \par
% We are also interested in \textit{irreducible, regular} monopoles on
% the four-dimensional pieces $X^\pm.$ Here we permit compact
% perturbations of the equations, which would readily rid us from the
% reducibles. Then again we restrict our attention to the open set of
% regular points.
%
% \par
% \ldots
\par
% The idea of the construction is as follows. 
Below is a quick survey of the gluing procedure.  
We put complete metrics on the two pieces $X^\pm,$ thus taking them
with infinite ends isometric to $[0,\infty)\times Y.$ There is a
limiting map which smoothly assigns a monopole on $Y$ to each
finite-energy monopole on the the \cyend \mfld. We cut each infinite
piece at place $\ell$ and glue the truncated manifolds $X^\pm_\ell$
along $Y.$ The resulting manifolds $X_\ell$ are diffeomorphic to $X$
and we try to successfully glue the monopoles on $X^\pm$ as the neck
is elongated.  We start naively by pasting together the monopoles
using a partition of unity.  Of course we need to add a correction
term to obtain an exact monopole; it would be the unique fixed point
of a contraction mapping on a Hilbert space. In fact, if
$\tilde{\xipt_\ell}$ is the approximate glued monopole and $\xihat$ is
the correction term, then we want the left-hand side of the expansion
$$\SW(\tilde{\xipt}_\ell+\xihat)=\SW(\tilde{\xipt}_\ell)+\D^1(\xihat)+Q(\xihat)$$
to be zero. Here $\D^1$ is the linearization of the $\SW$ map and
$Q$ is its second order approximation, both depending on
$\tilde{\xipt}_\ell,$ of course. The equation can be re-written as
$$\D^1(\xihat)=-\SW(\tilde{\xipt}_\ell)-Q(\xihat).$$
Thus if we can find a right inverse $\R$ for $\D^1$ and set
$\xihat=\R(\zetahat),$ then $\zetahat$ would be the fixed
point of the map
$$\F(\zetavector)=-\SW(\tilde{\xipt}_\ell)-Q(\R(\zetavector)).$$
\par
In the construction of a right inverse for $\D^1,$ we are led to
consider two complementary subspaces --- one being finite-dimensional
--- and work on each one separately. The main difficulty lies in
right-inverting $\D^1$ on the finite-dimensional subspace; this is
essentially due to the existence of obstruction spaces in the first
place. We will make estimates on the norms of certain operators,
obtained through Hodge theory, to eventually conclude that the desired
full right inverse can be constructed for $\ell$ sufficiently large.
The norm of this operator could grow exponentially in $\ell;$
nevertheless, the perturbation term will decline exponentially and
this conforms with the intuition that the approximate gluing is
increasingly ``better'' as $\ell$ becomes larger.
\par
The hard analysis culminates in the main theorem of this paper, whose
proof is completed in section \ref{sec:monopolesglue}. Use
$\M$ to denote the \sw moduli space 
% (the variants are distinguished by subscriprts) 
and let $\partial_\pm:\M(X^\pm)\to\M(Y)$ be the limiting map. 
Assume that $\M(Y)$ is Kuranishi-smooth and consider $\M^{irr}(Y),$ the
(smooth) irreducible part of the moduli space of $Y.$
Let $\M^*(X^\pm)$ consist of the regular points of the inverse image
of $\M^{irr}(Y)$ under the limiting map.
% Note that $\M^*(X^\pm)$ will not have reducible points either.
% Note that $\M^*(X^\pm)$ is already free from reducible points, too.
Note that this implies that $\M^*(X^\pm)$ is free from reducible
points, too. We continue to use the same notation for the restriction
of $\partial_\pm$ to $\M^*(X^\pm).$
\begin{thm*}
  Under the preceding assumptions on $X^\pm,$ $Y,$ and their moduli
  spaces, if the limiting maps
  $\partial_\pm:\M^*(X^\pm)\to\M^{irr}(Y)$ are transversal,
%   where $h^\pm$ are compactly-supported perturbation terms, 
  then there is an $L_0$ such that for each $\ell\geq 4L_0,$ the
  following holds. To any two regular monopoles $\xipt^+$ and
  $\xipt^-,$ respectively on $X^+$ and $X^-,$ with the same limiting
  value $\etapt\in\M^{irr}(Y),$ one can smoothly assign a monopole
  $\xipt_\ell$ on $X_\ell$ obtained through a canonical gluing scheme.
% This will produce a smooth gluing map 
% \ldots
\end{thm*}
%
% \par
% In the next section, we set our terminology and notation. The section
% afterwards is dedicated to the proof of the theorem above.
% \par
% \ldots
\par
This work is based on the author's PhD dissertation \cite{s}. 
\newpage
\section{Rudiments}
\label{sec:rudiments}%
% Stating preliminary facts and setting the notation and terminology in
% this section: moduli space, Chern-Simons, \ldots.
%
Let us first review some basic facts of \sw theory and meanwhile set
our notations along the way.
%
% \subsection{Seiberg-Witten Equations and Moduli}
% \label{sec:sw}%
\par
Let $X$ be a smooth, connected, oriented, riemannian \4mfld. We equip
$X$ with a \spinc-structure \s, 
\glossary{$\s$}%
that is a lifting of its principal tangent SO(4)-bundle $P$ 
\glossary{$P$}%
to a principal \spinc(4)-bundle $\tilde{P}.$
\glossary{$\tilde{P}$}%
Such liftings always exist and correspond (non-canonically) to classes
in $H^2(M,\bbZ)$ --- in fact, one can twist $\tilde{P}$ with any given
U(1)-bundle. Corresponding to Clifford representations of 
% \note{irreducible?}%
$\spinc(4)=\mbox{SU(2)}\times\mbox{SU(2)}\times\mbox{U(1)}/\{\pm 1\},$
we obtain the associated plus- and minus-spinor bundles $S^+$ and
\glossary{$S^+$}\glossary{$S^-$}
$S^-,$ which are 2-dimensional complex vector bundles with bundle
group U(2), as well as the determinant line bundle
$L=\det\tilde{P}=\det{S^+}=\det{S^-}.$ 
\glossary{$L$}%
\par
Any unitary connection $A$ 
\glossary{$A$}%
on the U(1)-bundle $L,$ in conjunction with the Levi-Civit\`{a}
connection on $X,$ will induce connections on the lifting $\tilde{P},$
as well as on $S^\pm.$ Thus a Dirac operator
$\dirac_A:\Gamma(S^\pm)\to\Gamma(S^\mp)$ can be defined by
\glossary{$\dirac_A$}\glossary{$\Gamma(V)$}%
$$\dirac_A(\Psi)=\sum_{j=1}^{n} e_j . \nabla_{e_j}(\Psi),$$
where $\{e_j\}_{j=1}^{n}$ is an orthonormal frame for $T_x X$ and
. denotes Clifford multiplication. The definition is frame-invariant.
\par
The \sw equations, in two unknowns $A$ and $\Psi,$ can now be written
as
$$
\left\{
  \begin{array}{l}
%   \label{eqn:sw}%
  F_A^+=\{\Psi\otimes\Psi^*\} \\
  \dirac_A(\Psi)=0,
  \end{array}
\right.
\eqno{(\SW)}
%\leqno{\mbox{(SW)}}
$$
where 
% $A$ is a unitary connection on $\L$ as mentioned above, 
\glossary{$\Psi$}%
$\Psi$ is a plus-spinor and the brackets denote the trace-less part of an
endomorphism. In other words, $\{\Psi\otimes\Psi^*\}$ denotes the
quadratic 
\glossary{$q$}%
$q(\Psi)=\Psi\otimes\Psi^*-\frac12|\Psi|^2\mbox{I}.$
In the same equation, $F_A^+$ denotes the self-dual part of the
curvature tensor under the Hodge $*$-operator. It defines a trace-free
representation of the
% even-grade, plus-
Clifford bundle $\Cl^+_0$ on $S^+$ via Clifford multiplication, thus
both sides of the first equation should be identified as trace-less
sections of the bundle of endomorphisms of plus-spinors, 
$\mend(S^+).$
\par
One can also consider the perturbed variant of the \sw equations
$$
\left\{
  \begin{array}{l}
%   \label{eqn:sw}%
  F_A^+=\{\Psi\otimes\Psi^*\}+ih \\
  \dirac_A(\Psi)=0,
  \end{array}
\right.
\eqno{(\SW_h)}
%\leqno{\mbox{(SW)}}
$$
where $h$ is a real self-dual 2-form on $X.$
\par
Sometimes it is convenient to consider the \sw \textit{map} on the
\textit{configuration space} $\C(X,\s)=\A(L)\times\Gamma(S^+)$
\glossary{$\C(X,\s)$}\glossary{$\A(V)$}%
\index{configuration space}%
consisting of a connection on the determinant line bundle and a
plus-spinor. The map
$\SW:\A(L)\times\Gamma(S^+)\rightarrow\forms2+(X)\times\Gamma(S^-)$ is
defined by
\glossary{$\SW$}%
$$\SW(A,\Psi)=(F^+_A-q(\Psi),\dirac_A(\Psi)).$$
$\SW_h$ is defined similarly.
\par 
We will feel free to make various assumptions on the configuration
spaces, for example by taking completions with respect to an
appropriate norm, or by considering only finite-energy configurations.
That should be clear from the context and we would invariably use the
same notation $\SW$ or $\SW_h.$
\par
The \textit{gauge group} $\G,$ 
\glossary{$\G(\tilde{P})$}%
\index{gauge group}%
i.e. the group of bundle-automorphisms of $\tilde{P},$ corresponds to
maps $X\to S^1.$ It right-acts on the configuration space $\C(X,\s),$
as well as on the solutions to the \sw equations $\S(X,\s),$ by
\glossary{$\S(X,\s)$}
$$(A,\Psi).g=(A+2g^{-1}dg, S^+(g^{-1})(\Psi)).$$ 
The stabilizer of $(A,\Psi)$ is trivial iff $\Psi\neq 0,$ in which
case the point is called {\it irreducible}. Reducible solutions have
the stabilizer $=S^1.$
\index{irreducible}\index{reducible}%
Dividing out the solution set by the action of the gauge group
produces the \sw \textit{moduli space} $\M(X,\s)=\S(X,\s)/\G.$
\glossary{$\M(X,\s)$}
\index{\sw moduli space}%
\par
A cohomological discussion of regularity is in order. To each point
$\xipt\in\C(X,\s)$ of the configuration space of a \4mfld $X,$ one can
assign the following diagram $\E|_\xipt$
\glossary{$\xipt$}\glossary{$\E|_\xipt$}%
% \index{elliptic complex}%
$$
0\longrightarrow
\forms03(X; i\bbR)
\stackrel{\D^0}{\longrightarrow}
\forms12(X;i\bbR)\oplus\Gamma_2(S^+)
\stackrel{\D^1}{\longrightarrow}
\Omega^2_{+,1}(X;i\bbR)\oplus\Gamma_1(S^-)
\longrightarrow 0,
\eqno{(\E|_\xipt)}
$$
\glossary{$\forms m k (Z;i\bbR)$}\glossary{$\Gamma_k(S^\pm)$}%
where $\forms m k (X;i\bbR)$ means the $L^2_k$-completion 
% \glossary{$L^2_k$}%
% \index{$L^2_k$-completion}%
(or a completion in another appropriate Sobolev norm, for that matter)
of purely imaginary $m$-forms with compact support and
$\Gamma_k(S^\pm)$ denotes the $L^2_k$-completion of the compactly
supported sections of the corresponding spinor bundles. $\D^0$ is the
linearization of the action of the gauge group and $\D^1$ is the
derivative of $\mbox{SW},$ both at the point $\xipt=(A,\Psi),$ i.e.
$$\D^0=(2d, -.\Psi);$$
\glossary{$\D^0$}%
$$
\D^1=
\left(
  \begin{array}{cc}
  d^+ & -Dq|_{\Psi} \\
  .\frac12\Psi & \dirac_{A}
  \end{array}
\right).
$$
\glossary{$\D^1$}%
Note that $\SW$ and $\SW_h$ are non-linear maps with the same
derivative $\D^1,$ so $\E|_\xipt$ remains unaltered with a
perturbation of the equations.
\par
Now if $\xipt$ happened to be a solution of $\SW$ or $\SW_h,$ then the
diagram $\E|_\xipt$ would be a complex; moreover it would even be an
\textit{elliptic} complex if $X$ were closed, or had appropriate
\index{elliptic complex}% \index{cohomology}%
boundary conditions, so it would have finite-dimensional cohomologies.
$H^0$ of such a complex turns out to be the tangent space to the
stabilizer
\index{stabilizer}%
of the gauge group action, $H^1$ is the Zariski tangent space
\index{Zariski tangent space}%
to the moduli space at $\xipt$ and $H^2$ is its obstruction space.
\index{obstruction space}%
By general Hodge theory,
\index{Hodge theory}%
these groups can be identified with the `harmonic forms`.
% \index{harmonic forms}
\par
Recall that a solution $\xipt=(A,\Psi)$ is called {\it irreducible}
\index{solution!irreducible}
if $\Psi$ is not identically zero; this is equivalent to $H^0(\E|_{\xipt})=0.$
% We say that the moduli space is {\it regular} (in an algebraic sense)
% at an irrdeucible point if the obstruction space at that point is
% trivial, i.e. $H^2(\E(Z)|_{\xipt})=0.$ Of course, according to
% Kuranishi picture, the moduli space is {\it smooth} (in a geometric
% sense) at such an irreducible point exactly when the Kuranishi map
% vanishes, even if the obstruction space is not trivial.
We call a solution %an irreducible point 
\textit{regular}
\index{solution!regular}%
(in an algebraic sense) if the obstruction space at that point is
trivial, i.e. $H^2(\E|_{\xipt})=0.$ 
Note that according to the Kuranishi picture,
% \index{Kuranishi picture}%
an irreducible solution is a {\it smooth} point
\index{solution!smooth}%
(in a geometric sense) if and only if the Kuranishi map
\index{Kuranishi map}%
vanishes, although the obstruction space may not be trivial.
\par
We will use superscripts to denote the terms of $\E.$ Depending on the
emphasis, when all else is clear from the context, we might use
combinations such as $\E(X),$ $\E(X,\s)|_\xipt$ and so on. Thus, for
example,
\begin{eqnarray*}
  \label{E}%
  \glossary{$\E^0(X)$}\glossary{$\E^1(X)$}\glossary{$\E^2(X)$}%
  \E^0(X)=\forms03(X; i\bbR) \quad , \quad
  \E^1(X)=\forms12(X;i\bbR)\plus\Gamma_2(S^+) \quad , \quad
  \E^2(X)=\Omega^2_{+,1}(X;i\bbR)\plus\Gamma_1(S^-) \quad
\end{eqnarray*}
We will also use boldface Greek letters (corresponding to the base
manifold) for points of the configuration space (for example $\xipt$
is a point of $\C(X)$) while the same Greek letters are used for the
vectors of the corresponding Zariski tangent spaces ($\xivector$
belongs to $H^1(\E|_\xipt)$). % $T^Z_\xipt\C$%
\par
It would be nice to review the setup for the case of a \textit{closed}
\4mfld $X.$ If we complete the configuration space and the gauge group
using the $L^2_2$ and $L^2_3$ Sobolev norms, respectively, then we
obtain an affine Hilbert space on which a Hilbert Lie group is acting.
% Thus we can form the quotient space $\M(X,\s)=\S(x,\s)/\G,$ which
% turns out to be Hausdorff. This is called the \sw \textit{moduli space}.
% Of course it 
The moduli space $\M(X,\s)$ happens to be Hausdorff, but it
might have singularities, for example when the action is
not free. It can be shown though that the solutions of a generic
perturbation of the \sw equations are all irreducible and regular,
%(more on regularity in section \ref{sec:complexes}), %
therefore resulting in a smooth moduli space $\M_h(X,\s).$ Using the
Atiyah-Singer index theorem and Bochner's formula, we conclude that
$\M_h$ is a finite-dimensional compact \mfld of formal dimension
$\frac14(c_1(L)^2-2\chi(M)-3\sigma(M)),$ where $\chi$ is the Euler
characteristic and $\sigma$ is the signature.
\par
The basic reference for the material so far is \cite{m}.
\par
One can mimic the preceding constructions on a \3mfld $Y$ with a
\spinc(3)-structure $\t$ to obtain analogs, 
% There are analogs for the preceding constructions on a \3mfld $Y$ with
% a \spinc(3)-structure, 
where there is only one spinor bundle $S$ and no self-duality. Thus,
for instance, one can define a map
$\SW^3:\A(L)\times\Gamma(S)\rightarrow\Omega^2(Y)\times\Gamma(S)$ by
\glossary{$\SW^3$}%
$$\SW^3(B,\Phi)=(F_B-q(\Phi),\dirac_B(\Phi)).$$
The moduli space for a closed \3mfld would generically be
of formal dimension zero, as the index of an elliptic operator on
an odd-dimensional closed manifold is zero.
%
% \par
% \ldots
\par
Another case of particular interest is when $X=\bbR \times Y$ is a
cylinder on a closed \3mfld $Y$ and $\xipt=\pi^*(\etapt)$ 
\glossary{$\etapt$}%
is a translation-invariant solution. Then the cohomologies of
$\E(X)|_{\xipt}$ can be identified in an obvious way with the
cohomologies of $\E(Y)|_{\etapt}$ below. (This is not the elliptic
complex that is officially associated to \3mflds --- it is
not even elliptic --- but it % officially=normally
would be more fitting to our discussion here.)
$$
\begin{array}{l}
\glossary{$\E(Y)|_{\etapt}$}%
0\cdr\forms03(Y; i\bbR)
\stackrel{\D^0}{\cdr}
\forms12(Y;i\bbR)\plus\Gamma_2(S)
\stackrel{\D^1}{\cdr}
\forms21(Y;i\bbR)\plus\Gamma_1(S)\cdr 0,
\end{array}
\eqno{(\E(Y)|_{\etapt})}
$$
where
$$\etapt=(B,\Phi);$$
$$\D^0=(2d, -.\Phi);$$
$$\D^1=
\left(
  \begin{array}{cc}
  d & -Dq|_{\Phi} \\
  .\frac12\Phi & \dirac_{B}
  \end{array}
\right) .
$$
%
%
%
% \par
% \ldots
\par
The situation is in general more complicated % would worsen
if $X$ is not closed, as we need controlling conditions near the
boundary or infinity, so as to keep the ellipticity. Therefore we work
with configurations with finite \textit{energy} when we deal with
manifolds with cylindrical ends.
%
% \bigskip
%
% \subsection{Chern-Simons-Dirac Functional}
% \label{sec:csd}%
So we continue to consider a cylinder $I\times Y,$ where $I=[c,d]$ is
an interval. The solutions to the \SW equations on this \4mfld turn
out to be the gradient flow lines
% \index{gradient flow line}%
of the so-called ``Chern-Simons-Dirac'' functional
\index{Chern-Simons-Dirac functional}%
on $\C(Y,\t)$:
$$
\mbox{CSD}(B,\Phi)=\int_Y F_{B_0}\wedge b +\frac{1}{2}\int_Y b\wedge db
+\int_Y \langle \Phi,\dirac_B\Phi\rangle d\mbox{vol},
\glossary{$CSD$}%
$$
where $b=B-B_0$ and $B_0$ is a fixed background connection \cite{mst}.
The singular points of this vector field, i.e. the static solutions,
correspond to solutions of $\SW^3.$ There are analogs for perturbed
equations, too. 
\par
The \textit{energy} of a solution $(A,\Psi)$
\index{energy of a solution}%
on a cylinder is defined by any of the following equivalent formulas.
We assume that $A$ is in temporal gauge
% \index{temporal gauge}
and we write $(A,\Psi)=(B(t),\Phi(t)),$ where $B(t)$ and $\Phi(t)$ are
connections and spinors on the \3mfld.
$$
\begin{array}{rcl}
\label{energy}
\glossary{$E$}
E(A,\Psi) & = & \int_I \|\dot{B}\|^2 +\|\dot{\Phi}\|^2 \medskip \\
%           & = & \int_I \|F_{B(t)}-q(\Phi(t))\|_Y^2+\|\dirac_{B(t)}(\Phi(t))\|_Y^2 \\
          & = & \int_I \|\nabla \mbox{CSD}(B(t),\Phi(t))\|^2  \medskip\\
          & = & \mbox{CSD}(B(d),\Phi(d))-\mbox{CSD}(B(c),\Phi(c)).
\end{array}
$$
\par
\index{end of a manifold}%
% by formal definition%
The \textit{end} of a \mfld is, formally, the inverse limit of its
co-compact subsets, ordered by inclusion. Intuitively, this is the
place where the manifold extends to infinity. We call a riemannian
\mfld $Z$ a \textit{\cyend manifold}
\index{\cyend \mfld}% 
\glossary{$Z$}%
if its end is orientation-preserving isometric to
$[0,\infty)\times Y,$ where $Y$ is an oriented, riemannian \3mfld.
\glossary{$Y$}%
We fix a smooth ``time coordinate'' function 
\glossary{$\tau$}%
$\tau\colon Z\to [-1,\infty)$ which agrees with the first coordinate
of $[0,\infty) \times Y$ on the end and is negative on the complement.
Given a positive real number $\ell>0,$ let
$Z_\ell=\tau^{-1}((-\infty,\ell]).$ 
\glossary{$Z_\ell$}%
For a pair of positive real numbers $0<\ell<\ell',$ let
$Z_{[\ell,\ell']}=\tau^{-1}([\ell,\ell']).$
\glossary{$Z_{[\ell,\ell']}$}%
\par
On \cyend manifolds, we will exclusively work with solutions with
finite energy on the ends for the sake of ellipticity.
% 
% $$
% \begin{array}{rcl}
% \label{energy}
% \glossary{$E$}
% E(A,\Psi) & = & \int_I \|\dot{B}\|^2 +\|\dot{\Phi}\|^2=
% \int_I \|F_{B(t)}-q(\Phi(t))\|_Y^2 +\|\dirac_{B(t)}(\Phi(t))\|_Y^2 \\
%           & = & \int_I \|\nabla f_{CSD}(\etapt(t))\|^2 \\
% % {\mbox{\qquad\qquad (where $\etapt(t)=(B(t),\Phi(t))$ is the corresponding flow line)}} \\
%           & = & f_{CSD}(\etapt(d))-f_{CSD}(\etapt(c)),
% %% \\
% %%             & = & \int_{I \times Y} F_A \wedge F_A .
% \end{array}
% $$
% where $\etapt(t)=(B(t),\Phi(t))$ is the corresponding flow line.
%
%%%
%
% To be written \ldots.
% %
% \subsection{Elliptic Complexes}
% \label{sec:complexes}%
% To be written \ldots.
% %
\newpage
\section{The Gluing Theorem}
\label{sec:gluing}%
\subsection{Gluing Cylindrical-End Manifolds}
\label{sec:manifoldglue}%
\glossary{$X^\pm$}\glossary{$\mend(X^\pm)$}\glossary{$X^o$}\glossary{$X^\#$}%
Let us fix two connected, oriented, \cyend riemannian \4mflds $X^\pm,$
each with a {\it single} end which is modeled on 
$[0,\infty)\times Y.$ We will also consider the cylinder
$X^o=\bbR\times Y.$ The notation $X^\#$ will then be used to denote
any of these three manifolds.
\par
\glossary{$\tau_\pm$}%
We will also introduce ``time'' coordinates $\tau$ on these manifolds
as follows. On $X^\pm,$ take the first coordinate map on the end
$[0,\infty)\times Y$ and choose any extension
$\tau_\pm:X^\pm\to[-1,\infty)$ which is identically -1 outside a
collar (the collar being identified with $(-1,0]\times Y$). On
$X^o=\bbR\times Y,$ $\tau_o$ is essentially the absolute value
function on the first coordinate, smoothed out at the origin. 
This will be made more precise further below.
% \newpage
%
\begin{figure}[h]
  \refstepcounter{figure}       % Fixes the pre-cross-referencing problem.
  \label{fig:cyend}%
  \addtocounter{figure}{-1}     % Fixes the pre-cross-referencing problem.
  \begin{center}
    \leavevmode
      \ForceWidth{11cm}
      \BoxedEPSF{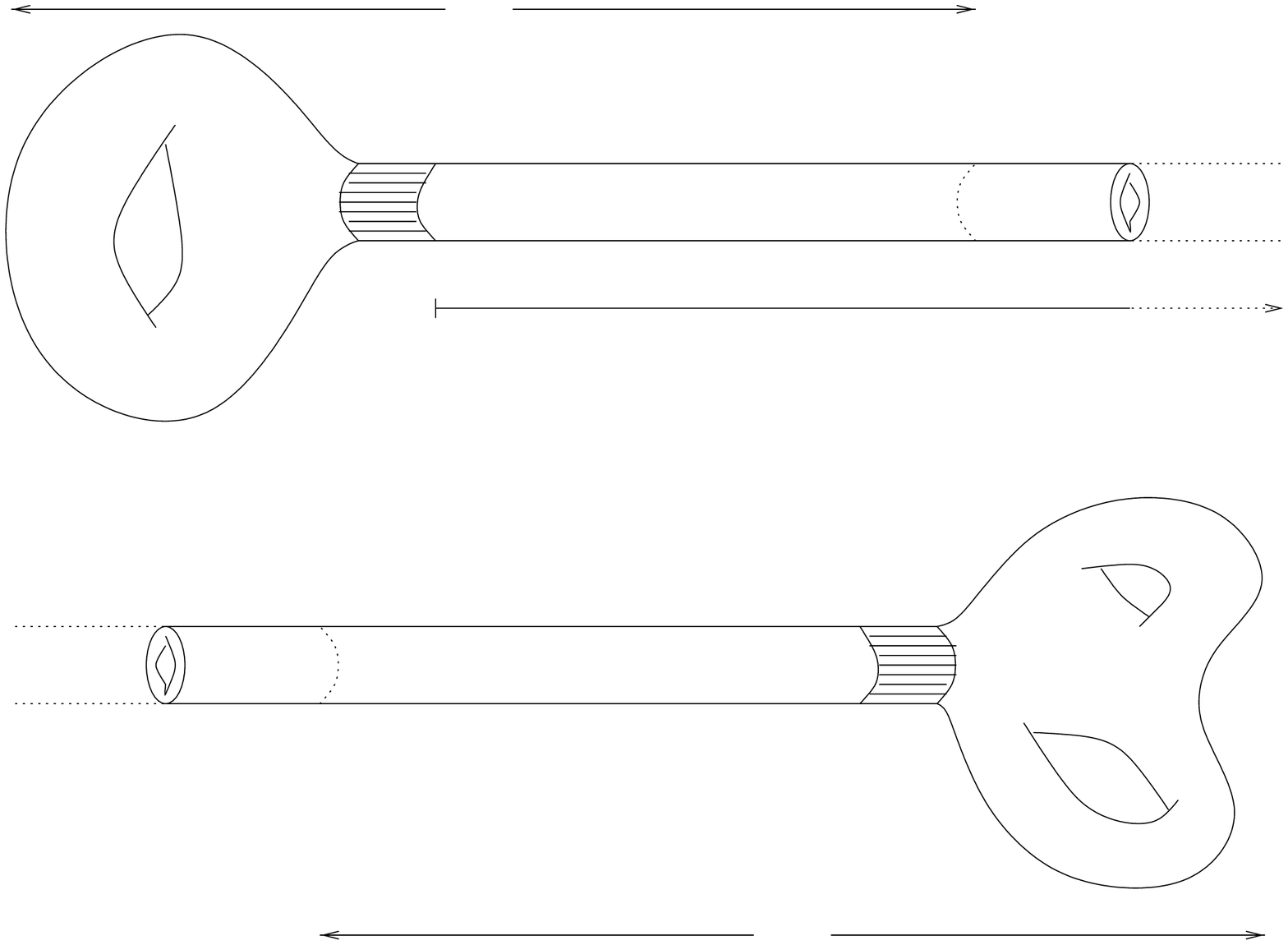}
       \begin{picture}(0,0)(395,367)
       {\tiny
        \put(189,481){$X^+_\ell$} \put(182,417){0} \put(300,417){$\tau_+=\ell$}
        \put(360,434){$Y$} \put(240,414){$\mend(X^+)$}
        \put(103,322){$Y$} \put(145,304){$\tau_-=\ell$} \put(285,304){0}
        \put(265,255){$X^-_\ell$}
       }
       \end{picture}
    \caption{Two \cyend manifolds $X^+_\ell$ and $X^-_\ell$}
  \end{center}
\end{figure}
\glossary{$X^\pm$}\glossary{$X^o$}\glossary{$\tau_\pm$}\glossary{$X^\pm_\ell$}%
% \bigskip % \vspace{1cm}
%
\begin{figure}[h]
  \refstepcounter{figure}       % Fixes the pre-cross-referencing problem.
  \label{fig:glued}%
  \addtocounter{figure}{-1}     % Fixes the pre-cross-referencing problem.
  \glossary{$X^\pm_\ell$} \glossary{$X_\ell$} \glossary{$L_0$}%
  \glossary{$C_{\ell'}$} \glossary{$\tau_o$} \glossary{$\tau_\ell$}%
  \begin{center}
    \leavevmode
      \ForceWidth{15cm}
      \BoxedEPSF{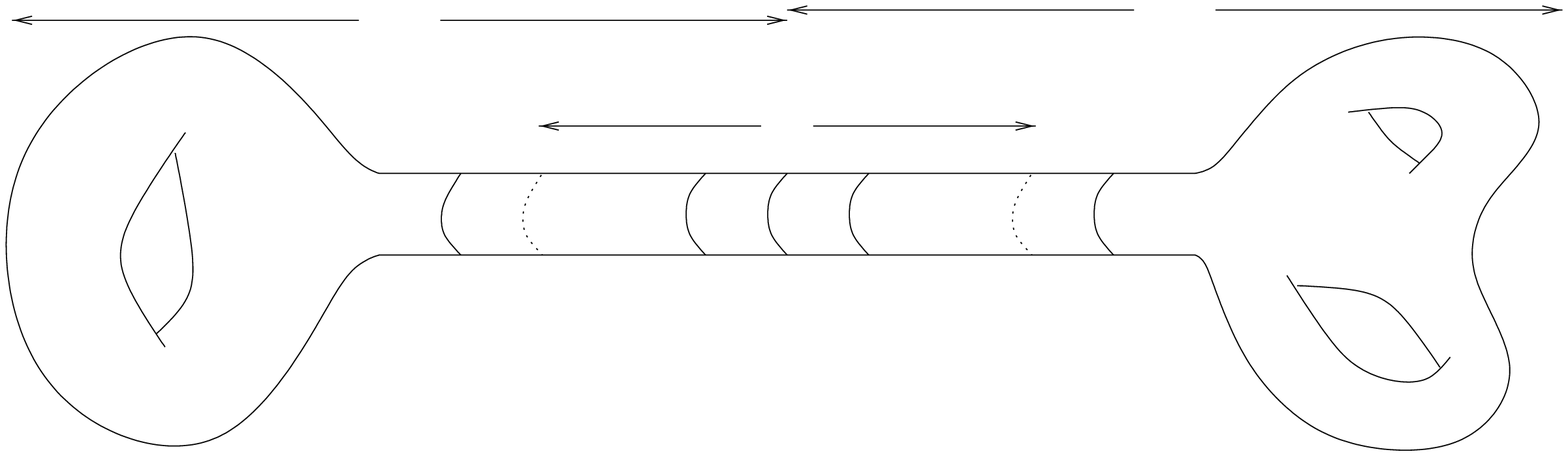}
       \begin{picture}(0,0)(235,59) %(452,118)
       {\tiny
        \put(125,177){$X^+_\ell$} \put(338,180){$X^-_\ell$}
        \put(231,148){$C_{\ell'}$} \put(232,105){0}
        \put(253,105){$L_0$} \put(200,105){$-L_0$}
        \put(299,105){$\ell'$} \put(162,105){$-\ell'$}
        \put(320,105){$\ell$} \put(141,105){$-\ell$}
       }
       \end{picture}
    \caption{$X^+_\ell$ and $X^-_\ell$ glued together to form $X_\ell$}
  \end{center}
\end{figure}
\bigskip
\par
\glossary{$X_\ell$}\glossary{$X^\pm_\ell$}\glossary{$X$}%
\glossary{$\ell$}\glossary{$L_0$}\glossary{$C_{\ell'}$}%
We now form a family of \4mflds $X_\ell,$ $\ell>0,$ as follows. First
truncate the manifolds $X^\pm$ at $\tau_\pm=\ell$  % $(\tau_\pm)^{-1}(\ell)$
to produce $X^\pm_\ell$ consisting of points with
$\tau_\pm(x)\leq\ell.$ We then obtain $X_\ell$ by gluing $X^+_\ell$
and $X^-_\ell$ along their boundaries (see figures~\ref{fig:cyend} 
and \ref{fig:glued}). We are interested in $X_\ell$ for $\ell$ large % $\ell\gg 1$ 
and we would eventually assume $\ell>4L_0,$ where $L_0$ is a large,
but fixed, positive number. The manifolds $X_\ell$ are just
diffeomorphic versions of one and the same manifold $X,$ being
elongated along a tube. We will also re-parameterize the long cylinder
$C_\ell=\mend(X^+_\ell)\cup \mend(X^-_\ell)$ inside $X_\ell,$
identifying it with $[-\ell,\ell]\times Y$ as in
figure~\ref{fig:glued}. $C_{\ell'}$ is then the chunk of $C_\ell$
parameterized as $[-\ell',\ell']\times Y.$
\par
\newenvironment{tau1}{}{}%
\begin{tau1}\label{tau-page}\end{tau1}%
\glossary{$\tau_\ell$}\glossary{$\tau_o$}%
Now we get back to our discussion of the time coordinate and, in the
mean time, we also introduce a time-coordinate function
% \index{time-coordinate function}
$\tau_\ell$ on the manifold $X_\ell.$ It is identical to $\tau_+$ on
$X^+_{\ell-1}$ and to $\tau_-$ on $X^-_{\ell-1}$, and smoothly
interpolates between the two such that its value on $C_1$ is in the
interval $[\ell-1,\ell].$ This choice of interpolation can be made
independently of $\ell$ and therefore the function $\ell-\tau_\ell$
converges, uniformly on compact subsets of $X^o,$ to the function
$\tau_o$ alluded to earlier. Outside of $C_1=[-1,1]\times Y$ we have
$\tau_o(t,y)=|t|.$ 
% for $(t,y)\in X^o=\bbR\times Y.$
% \note{For $\tau_\ell$ \& $\tau_o,$ see [MM], pp.  9-10. There,
% $\tau_{cyl}=$ our $\tau_o.$}
\par
We are also able to patch \spinc-structures $\s^\pm$ (on $X^\pm$)
which agree on the ends. More precisely, there is a principal
\spinc(3)-lifting $\tilde{Q}$ of the principal orthonormal tangent
bundle of $Y$ such that $\tilde{P}^\pm|_{\mend(X^\pm)} \iso
\bbR\times\pi^*\tilde{Q},$ where $\pi:I\times Y \rightarrow Y$ is the
projection. Indeed, on $\mend(X^\pm),$ there is an embedding
$\pi^*\tilde{Q} \hookrightarrow \tilde{P}^\pm,$ induced from the lift
of the obvious embedding $SO(3) \hookrightarrow SO(4).$ See \cite{mst}
for details.  By fixing the isomorphism above, we can form a
\spinc-structure $\s_\ell$ on $X_\ell$ which is compatible with the
original \spinc-structures.
\subsection{Approximate Gluing of Monopoles}
\label{sec:approxglue}%
\glossary{$\xipt^\pm$}%
% \note{!Care!}%
% \note{Can we instead talk of $\SW_h$ from the beginning?}%
Let $\xipt^\pm=(A^\pm,\Psi^\pm)$ be finite-energy solutions to the %perturbed
\sw equations $\SW$ on $X^\pm.$ 
Based on our earlier assumptions, $\xipt^\pm$ will be regular and
irreducible.
% Let us choose $h_\pm$ generically so that the moduli space of
% solutions is irreducible [and smooth??] at each and every point.
% \note{correct the assumptions}
\par
\glossary{$\kappa$}%
By a result of Simon \cite{mmr}, ``finite energy implies finite
length'' for the solution, which is now viewed as a gradient flow line
for $CSD$ on the cylindrical end. It further implies ``exponential
decay'' to a solution $(B,\Phi)$ of $\SW^3$ on $Y.$
% $(A_\infty, \Psi_\infty)$ 
See \cite{mmr}. The exponent $\kappa$ in this exponential decay is
less than half the minimum of the absolute value of the eigenvalues of
the Hessian
% \note{for $\kappa,$ see \cite{mm}, pp.24, 49}
$$
\mbox{Hess}(CSD)=
\left(
  \begin{array}{cc}
    *d & Dq|_{\Psi} \\
    0  & \dirac_A
  \end{array}
\right),
$$
which is the linearization of the gradient flow at a critical point.
Therefore $\kappa$ has a bound which simply depends on the
eigenvalues of $\Delta$ and $\dirac_A.$
\par
We also assume that $\xipt^\pm$ converge to the same irreducible
solution $\etapt$ on $Y.$
% \note{$\etapt$ reducible, smooth?}%
This is tantamount to considering $(\xipt^+,\xipt^-)$ as a point in
the fiber product 
\index{fiber product}%
$\M^*(X^+)\times_U\M^*(X^-),$ defined as the pull-back of the diagram
$$
\xymatrix
{
{\M^*(X^+)}\times_U{\M^*(X^-)} \ar@{.>}[r] \ar@{.>}[d] &
{\M^*(X^-)|_U} \ar[d]^{\partial_-} \\
{\M^*(X^+)|_U} \ar[r]^{\partial_+} & **[r] U \subset \M^{irr}(Y),
}
$$
\glossary{$\partial_\pm$}%
where $\partial_\pm:\M^*(X^\pm)\rightarrow\M^{irr}(Y)$ are the
\index{limiting map}%
limiting maps for the flow lines, $U$ an arbitrary (smooth) \nbhd of
% \note{$U$ smooth?}%
$\etapt$ in $\M^{irr}(Y)$ and $\M^*(X^\pm)|_U=(\partial_\pm)^{-1}(U).$
$\M$'s may denote any variant of the \sw moduli spaces, depending on
the context.
\par
We will later also need to assume a transversality condition at
$\etapt.$ % $\xipt_\infty.$ 
%
% \vspace{2cm}
\par
\medskip
\glossary{$\xipt^\pm_\ell$}\glossary{$\xipt_\ell$}%
\glossary{$\tilde{\xipt}_\ell$}\glossary{$\lambda$}%
\glossary{$\tilde{\gamma}=\tilde{\gamma}_\ell$}%
While passing from $X^\pm$ to $X^\pm_\ell,$ we truncate $\xipt^\pm$ to
$\xipt^\pm_\ell$ as well. Then we glue $X_\ell^+$ and $X_\ell^-$
together to form $X_\ell$ and our goal is to ``glue'' $\xipt^+_\ell$
and $\xipt^-_\ell$ to construct a solution $\xipt_\ell$ to
$\SW$ on the glued-up manifold $X_\ell,$ for large $\ell.$ 
% Here $h=h_++h_-.$
%
To this end, our first step would be to construct an \textit{approximate} solution
$\tilde{\xipt}_\ell=(\tilde{A}_\ell,\tilde{\Psi}_\ell)$ of $\SW$ on
$X_\ell,$ using a partition of unity $\{\lambda, 1-~\lambda\}$ which is
constant outside of $C_{2L_0}.$ Thus we can define an ``approximate
gluing map''
\index{gluing map!approximate}%
$$\tilde{\gamma}=\tilde{\gamma}_\ell:
\S(X^+,\s^+)\times\S(X^-,\s^-)\to\C(X_\ell,\s_\ell)$$
by
$\tilde{\xipt}_\ell:=\tilde{\gamma}_\ell(\xipt^+,\xipt^-)=
\lambda\xipt^+_\ell+(1-\lambda)\xipt^-_\ell.$
\par
\medskip 
\glossary{$L^2_{k,-\delta}$}\glossary{$\|\alpha\|_{k,-\delta}$}%
\index{weighted Sobolev norm}%
Due to technical problems arising from the presence of an obstruction
space $H^2,$ we shall use \textit{weighted} Sobolev norms to allow
some small exponential growth on forms and spinors.  Let $\alpha$ be a
$C^\infty,$ compact supported $m$-form on a \cyend \mfld $Z,$ $\tau$
the time coordinate on $Z,$ $\nabla$ the Levi-Civit\`{a} connection
and $\delta$ a real number. Define the $L^2_{k,-\delta}$-norm of
$\alpha$ as
$$\|\alpha\|_{k,-\delta}=
\left(\sum_{j=0}^{k}\int_Ze^{-\delta\tau}|\nabla^j\alpha|^2\right)^{1/2}.$$
\index{Levi-Civit\`{a} connection}%
We will denote the $L^2_{k,-\delta}$-completion of the space of
$C^\infty,$ compact supported $m$-forms on $Z$ by
$\wforms{m}{k}{-\delta}(Z).$ Analogous terminology will also be used
for spinor fields, except that a hermitian connection on the spinor
bundle must be used instead of the Levi-Civit\`{a} connection.
% \par
On $X_\ell,$ we can define the weighted norms similarly; $\tau_\ell$
is to be used instead of $\tau.$ However, note that Sobolev norms
with different weights are equivalent on a closed manifold.
\par
We have the following estimate.
\begin{prop}
  \label{exp-decay}%
  \glossary{$\tilde{C}$}%
  For any $\delta\geq 0$ and any $\ell>4L_0,$
  $$\|\SW(\tilde{\xipt}_\ell)\|_{1,-\delta}\leq
  \tilde{C}e^{-(\kappa+\frac{\delta}{2})(\ell-2L_0)},$$
% \note{\cite{mm}, p.27.}
  for some constant $\tilde{C}$ which is independent of $\ell$ and
  $L_0.$ % \endpf
\end{prop}
\begin{proof} % {\thetheorem}
  $\SW(\tilde{\xipt}_\ell)$ is zero when $\lambda$ is constant.
  So we only need to estimate it on $C_{2L_0}.$ Using the fact that
  $\SW(\xipt^+)$ and $\SW(\xipt^-)$ are both zero, we get
%   we find out that $\SW_{h}(\tilde{\xipt}_\ell)$ is actually
%   independent of $h$ and we have
  $$\SW(\tilde{\xipt}_\ell)= (-d^+(\lambda a)-
  \lambda^2q(\psi)-Dq(\Psi^-,\lambda\psi), -\nabla(\lambda)\psi),$$
  where $a=A^+-A^-$ and $\psi=\Psi^+-\Psi^-$ on $C_{2L_0}.$ An easy
  computation shows that $|q(\psi)|=\frac12|\psi|^2$ and that
  $|Dq(\psi,\psi')|\leq2|\psi||\psi'|.$ Now the fact that the
  solutions decay exponentially fast with exponent $\kappa$ gives the
  desired estimate.
\end{proof}
\par
\glossary{$\R$}%
Next, we deform $\tilde{\xipt}_\ell$ to a solution of
$\SW(\xipt_\ell)=0.$ For this, as was pointed out in
section \ref{sec:intro}, we will need a right inverse $\R$
for $\D^1$ in the following diagram at the point
$\tilde{\xipt}_\ell=(\tilde{A}_\ell,\tilde{\Psi}_\ell).$
$$
0\cdr\wforms03{-\delta}(X_\ell; i\bbR)
\stackrel{\D^0}{\cdr}
\wforms12{-\delta}(X_\ell;i\bbR)\plus\Gamma_{2,-\delta}(S^+)
\stackrel{\D^1}{\cdr}
\Omega^2_{+,1,-\delta}(X_\ell;i\bbR)\plus\Gamma_{1,-\delta}(S^-)\cdr 0
$$
\par
This diagram is \textit{not} even a complex, since
$\dirac_{\tilde{A}_\ell}(\tilde{\Psi}_\ell) \neq 0.$
\par
We will construct such an inverse using the chain homotopies that
already exist between the complexes of each piece of our manifold and
their cohomologies. Let $0<\delta<\kappa/2$ and consider the following
diagram of two copies of the complex $\E_{-\delta}(X^\pm)$ at
$\xipt^\pm.$ 
\glossary{$\E_{-\delta}(X^\pm)$ at $\xipt^\pm$}%
$$
\glossary{$=\E_{-\delta}(X^\pm)|_{\xipt^\pm}=\E_{-\delta}(X^\pm)$}%
\xymatrix
{0 \ar[r] & \E^0_{-\delta}(X^\pm) \ar[r]^{\D^0} &
   \E^1_{-\delta}(X^\pm)
   \ar[r]^{\D^1} \ar[dl]^(0.37){\L_\pm}
   \ar[d]^{\Pi^1_\pm} &
%   \ar[r]^{\D^1} \ar[dl]^(0.33){{\displaystyle \L_\pm}}
%   \ar[d]^{{\displaystyle \Pi^1_\pm}} &
   \E^2_{-\delta}(X^\pm)
   \ar[r] \ar[dl]^(0.37){\R_\pm} & 0 \\
%   \ar[r] \ar[dl]^(0.33){{\displaystyle \R_\pm}} & 0 \\
 0 \ar[r] & \E^0_{-\delta}(X^\pm) \ar[r]^{\D^0} &
   \E^1_{-\delta}(X^\pm)
   \ar[r]^{\D^1} &
   \E^2_{-\delta}(X^\pm)
   \ar[r] & 0,
}
$$
where
$$
\begin{array}{rcl}
  \E^0_{-\delta}(X^\pm)&=&\wforms03{-\delta}(X^\pm; i\bbR),\\
  \E^1_{-\delta}(X^\pm)&=&\wforms12{-\delta}(X^\pm;i\bbR)\plus\Gamma_{2,-\delta}(S^+), \\
  \E^2_{-\delta}(X^\pm)&=&\Omega^2_{+,1,-\delta}(X^\pm;i\bbR)\plus\Gamma_{1,-\delta}(S^-),
\end{array}
$$
using weighted Sobolev completions.
% \par
These are the complexes associated to solutions $\xipt^\pm$ of
$\SW$ on $X^\pm.$ 
% The perturbation terms $h_\pm$ are chosen suitably so that
% \note{!Care!}%
As $\xipt^\pm$ are regular and irreducible, 
the complexes above have no cohomology except possibly in degree one.
$\Pi^1_\pm$ is the projection onto this cohomology, represented by the
harmonic 'one-forms', and the parametrices $\L_\pm$ and $\R_\pm$ are
constructed using Hodge theory,
\index{Hodge theory}%
so that 
% in a way that the following relations hold between the operators.
$$
\glossary{$\Pi^1_\pm$}\glossary{$\L_\pm$}\glossary{$\R_\pm$}%
\left\{
  \begin{array}{l}
  \L_\pm\circ\D^0={\mathrm I} \\
  \D^0\circ\L_\pm+\R_\pm\circ\D^1={\mathrm I}-\Pi^1_\pm \\
  \D^1\circ\R_\pm={\mathrm I}.
  \end{array}
\right.
$$
% \glossary{$\R_+$}\glossary{$\R_-$}\glossary{$\R_o$}%
Note that $\R_\pm$ are right inverses for $\D^1$ on the respective
pieces.
\par
Now recall that the solutions $\xipt^\pm,$ considered now as
gradient flow lines, both converge to the same static solution
$\etapt.$ Thus, similarly as above, consider the corresponding complex
% \glossary{$\E(Y)$ at $\etapt$}
\glossary{$\E_{\delta}(X^o)$ at $\xipt^o$}%
$\E_{\delta}(X^o)$ %$=\E_{-\delta}^{X^o}(\xipt^o)$
on the cylinder $X^o=\bbR\times Y$ at the constant solution 
% $\xipt^o$ 
(i.e. the pull-back of $\etapt$ to $X^o$), still denoted by $\etapt.$
% and $0<\delta<\kappa/2.$ 
% \note{$\etapt$ should be irreducible here!}%
$$
\glossary{$\E_{\delta}(X^o)|_{\xipt^o}$}%
\xymatrix
{0 \ar[r] & \E^0_{\delta}(X^o) \ar[r]^{\D^0} &
   \E^1_{\delta}(X^o)
   \ar[r]^{\D^1} \ar[dl]^(0.37){\L_o}
   \ar[d]^{\Pi^1_o} &
%   \ar[r]^{\D^1} \ar[dl]^(0.33){{\displaystyle \L_o}}
%   \ar[d]^{{\displaystyle \Pi^1_o}} &
   \E^2_{\delta}(X^o)
   \ar[r] \ar[dl]^(0.37){\R_o}
   \ar[d]^{\Pi^2_o} & 0 \\
%   \ar[r] \ar[dl]^(0.33){{\displaystyle \R_o}}
%   \ar[d]^{{\displaystyle \Pi^2_o}} & 0 \\
 0 \ar[r] & \E^0_{\delta}(X^o) \ar[r]^{\D^0} &
   \E^1_{\delta}(X^o)
   \ar[r]^{\D^1} &
   \E^2_{\delta}(X^o)
   \ar[r] & 0,
}
$$
where, similarly,
$$
\begin{array}{rcl}
  \E^0_{\delta}(X^o)&=&\wforms03{\delta}(X^o; i\bbR),\\
  \E^1_{\delta}(X^o)&=&\wforms12{\delta}(X^o;i\bbR)\plus\Gamma_{2,\delta}(S^+), \\
  \E^2_{\delta}(X^o)&=&\Omega^2_{+,1,\delta}(X^o;i\bbR)\plus\Gamma_{1,\delta}(S^-),
\end{array}
$$
and,
$$
\left\{
  \begin{array}{l}
  \L_o\circ\D^0={\mathrm I} \\
  \D^0\circ\L_o+\R_o\circ\D^1={\mathrm I}-\Pi^1_o \\
  \D^1\circ\R_o={\mathrm I}-\Pi^2_o.
  \end{array}
\right.
$$
Note here that $\Pi^2_o$ is an obstruction for $\R_o$ to be a right
inverse to $\D^1.$
\par
\glossary{$\tilde{\L}$}\glossary{$\tilde{\R}$}\glossary{$\tilde{\Pi}^1$}\glossary{$\tilde{\Pi}^2.$}%
\glossary{$\mu_+$}\glossary{$\mu_-$}\glossary{$\mu_o$}%
Once again, we will be using partitions of unity to splice these
parametrices and projections together. We will pick a partition 
% \note{\cite{mm}, p. 59,}%
$\{\mu_+^2,\mu_o^2,\mu_-^2 \}$ on $X_\ell$ such that each $\mu$ satisfies
$|\nabla^n\mu| \leq {(\frac2{L_0})}^n.$ $\mu_+$ is supported on
$X^+_\ell$ and is constant outside $X_{[-2L_0,-L_0]}.$
(See Figure 3.) Symmetrically, $\mu_-$ is supported on $X^-_\ell$ and is
constant outside $X_{[L_0,2L_0]}.$ $\mu_o$ is therefore supported on $C_{2L_0}.$
\setcounter{figure}{2}
\begin{figure}[h]
  \label{fig:partition}%
  \glossary{$\mu_+$}\glossary{$\mu_-$}\glossary{$\m,u_o$}%
  \begin{center}
    \leavevmode
       \ForceWidth{15cm}
       \BoxedEPSF{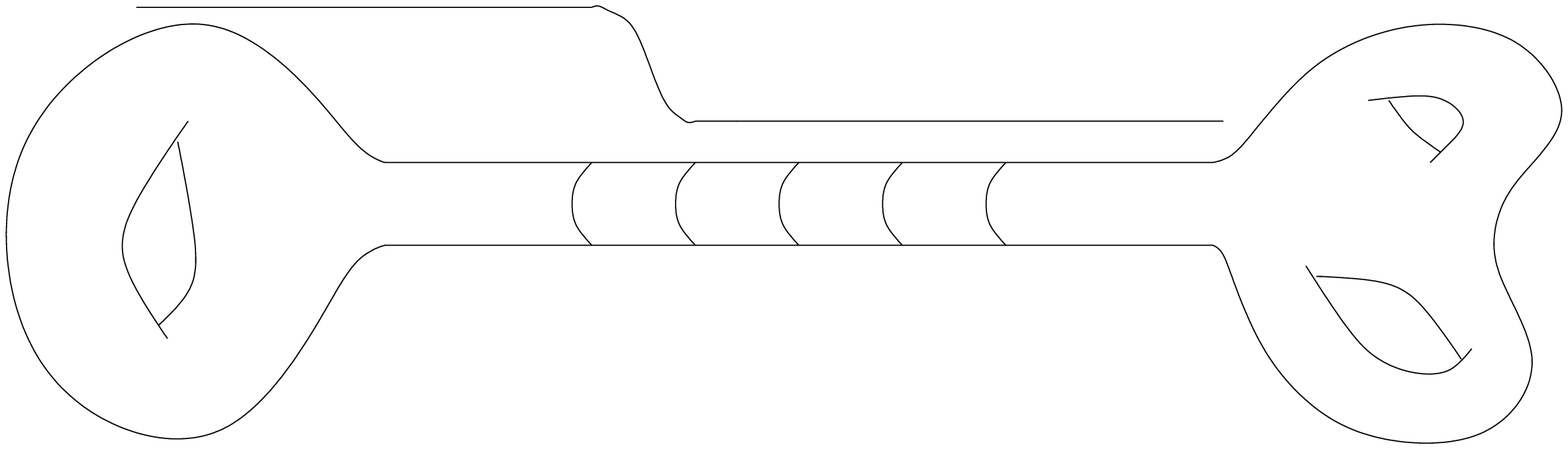}
        \begin{picture}(0,0)(232,55)
        {\tiny
        \put(130,180){1} \put(235,150){0} \put(195,170){$\mu_+$}
        \put(30,120){$X^+_\ell$} \put(395,120){$X^-_\ell$}
        \put(235,100){0} \put(260,100){$L_0$} \put(287,100){$2L_0$}
        \put(200,100){$-L_0$} \put(165,100){$-2L_0$}
        }
        \end{picture}
    \caption{Graph of $\mu_+.$ That of $\mu_-$ (not drawn) is the
    mirror image of $\mu_+$ on the right.}
  \end{center}
\end{figure}
\\
\glossary{$\tilde{\R}$}\glossary{$\tilde{\L}$}\glossary{$\tilde{\Pi}^1$}\glossary{$\tilde{\Pi}^2$}%
% Now we define an approximate right inverse
% \par
Now we paste the parametrices $\R_\#$'s, $\L_\#$'s and the projections
$\Pi_\#$'s on the pieces, using this partition of unity, to produce
$\tilde{\R},$ $\tilde{\L},$ $\tilde{\Pi}^1$ and $\tilde{\Pi}^2.$ As
we will be more interested in $\tilde{\R}$ and $\tilde{\Pi}^2,$
% $\tilde{\R}=\tilde{\R}_\ell(\xipt^+,\xipt^-)$ 
% and
% $\tilde{\Pi}^2=\tilde{\Pi}^2_\ell(\xipt^+,\xipt^-),$
we will give their explicit definitions below --- of course
$\tilde{\L}$ and  $\tilde{\Pi}^1$ are defined in a similar way.
\glossary{$\zetavector$}%
For  $\zetavector\in\E^2_{-\delta}(X_{\ell}),$
$$
\begin{array}{l}
  \label{eq:approx-ops}%
  \tilde{\R}\zetavector=\mu_+\R_+(\mu_+\zetavector)+\mu_o\R_o(\mu_o\zetavector)+\mu_-\R_-(\mu_-\zetavector), \\
  \tilde{\Pi}^2\zetavector= \mu_o \Pi^2_o(\mu_o\zetavector).
\end{array}
\addtocounter{equation}{1}
\eqno{(\theequation)}
$$
\par
\index{right inverse!approximate}%
These glued operators \textit{approximately} serve as their
counterparts on each piece, in the sense of the following lemma.
$\tilde{\R}=\tilde{\R}_\ell(\xipt^+,\xipt^-):\E^2_{-\delta}(X_\ell)
\to\E^1_{-\delta}(X_\ell)$ is supposed % going 
to be our first approximation of a right inverse for
$\D^1:\E^1_{-\delta}(X_\ell)\to\E^2_{-\delta}(X_\ell).$ 
% \note{$\E_{-\delta}(X_\ell)$?}%
% \note{define $\tilde{\Pi}^2$ explicitly.}%
%
\begin{lemma}
\label{first-estimates}%
\glossary{$\tilde{K}$}%
There is a constant $\tilde{K}$ such that
$$
\left\{
  \begin{array}{l}
    \|\tilde{\L}\circ\D^0-{\mathrm I}\| \leq \frac{\tilde{K}}{L_0} \\
    \|\D^0\circ\tilde{\L}+\tilde{\R}\circ\D^1-{\mathrm I}+\tilde{\Pi}^1\| \leq \frac{\tilde{K}}{L_0} \\
    \|\D^1\circ\tilde{\R}-{\mathrm I}+\tilde{\Pi}^2\| \leq \frac{\tilde{K}}{L_0}.
  \end{array}
\right.
$$ % \endpf
\end{lemma}
\begin{proof}
\glossary{K',K''}\glossary{$\wedgedot$}%
We prove the last estimate; the others are proved similarly. Using
the fact that $\D^1\circ\R_\pm={\mathrm I}$ and $\D^1\circ\R_o={\rm  I}-\Pi^2_o,$ 
we can write, using operator commutators,
$$
\begin{array}[t]{lll}
  \D^1\circ\tilde{\R}(\zetavector) &=&
  [\D^1,\mu_+]\R_+(\mu_+\zetavector)+
  [\D^1,\mu_o]\R_o(\mu_o\zetavector)+
  [\D^1,\mu_-]\R_-(\mu_-\zetavector) \\ &&
  +\mu_+^2\zetavector+\mu_o^2\zetavector+\mu_-^2\zetavector
  -\mu_o\Pi^2_o\mu_o\zetavector
\end{array}
$$
for $\zetavector\in\E^2_{-\delta}(X_{\ell}).$ 
% This is because $\D^1\circ\R_o={\mathrm I}-\Pi^2_o.$ 
Therefore, 
$$\D^1\circ\tilde{\R}(\zetavector)-\zetavector+\tilde{\Pi}^2(\zetavector)=
  [\D^1,\mu_+]\R_+(\mu_+\zetavector)+
  [\D^1,\mu_o]\R_o(\mu_o\zetavector)+
  [\D^1,\mu_-]\R_-(\mu_-\zetavector).$$
Thus to estimate $\|\D^1\circ\tilde{\R}-{\mathrm I}+\tilde{\Pi}^2\|$ we
should estimate the commutator norms. Note that
$\D^1(f\xivector)=f\D^1(\xivector)+df\wedgedot\xivector,$ where $f$
is a scalar function and for a one-form $\omega$ on $X_\ell,$
$\omega\wedgedot(a,\psi)=((\omega\wedge a)^+,\omega\,.\,\psi),$ where
in the second component dot denotes Clifford multiplication.  As a
result, the commutator $[\D^1,\mu_\#]\R_\#(\mu_\#\zetavector)$ is just
$d\mu_\#\wedgedot\R_\#(\mu_\#\zetavector)$ and, using
$|\nabla\mu|\leq\frac2{L_0},$ we obtain
$$
\begin{array}[t]{lll}
  \|[\D^1,\mu_\#]\R_\#(\mu_\#\zetavector)\| &=&
  \|d\mu_\#\wedgedot\R_\#(\mu_\#\zetavector)\| \\
  &\leq&
  \frac{\tilde{K}'}{L_0}\|\R_\#(\mu_\#\zetavector)\| \\
  &\leq&
  \frac{\tilde{K}'}{L_0}\|\R_\#\|.\|\mu_\#\zetavector\| \\
  &\leq&
  \frac{\tilde{K}''}{L_0}\|\R_\#\|.\|\zetavector\|,
\end{array}
$$
where we have everywhere used the $L^2_{1,-\delta}$ weighted Sobolev
norm. The last inequality is justified by the Sobolev embedding
$L^2_{3,0}\otimes L^2_{1,-\delta}\hookrightarrow L^2_{1,-\delta}.$
Now take $\tilde{K}=\frac{\tilde{K}''}{L_0}\|\R_\#\|.$
\end{proof}
\par
Unfortunately, $\tilde{\Pi}^2$ is \textit{not} a projection; however,
it gets closer to one as $L_0 \to \infty$ and it has a right inverse.
\begin{lemma}
\label{tildePi2}%
\glossary{$\tilde{\Pi}^2$}\glossary{$K_2$}%
  The operator $\tilde{\Pi}^2$ satisfies
  $\|\tilde{\Pi}^2-\tilde{\Pi}^2\circ\tilde{\Pi}^2\|_{1,-\delta} \leq
  K_2 e^{-\delta L_0}$ and has a right-inverse, defined on
  $\im(\tilde{\Pi}^2),$ whose operator norm is at most $K_2.$
  % \endpf
\end{lemma}
\begin{proof}
\glossary{$\wedgedot$}\glossary{\beta}%
  We are going to calculate $\tilde{\Pi}^2\circ\tilde{\Pi}^2.$ For
  this, we first find an expression for $\Pi^2_o,$ which is the
  projection $\E^2_{\delta}(\bbR\times Y)\to\H^2_{\delta}(\bbR\times
  Y)$ onto the harmonic forms. Let $h_1,\cdots,h_n$ be an orthonormal
  base for $\H^1(\E(Y)).$ Using our notation $\wedgedot$ from the
  previous lemma, we construct an isometry
%   There is an isometry
  $\H^1(\E(Y))\mathop{\longrightarrow}\limits^{\simeq}
  \H^2_{\delta}(\bbR\times Y)$ given by $\etavector\mapsto
  ce^{-\delta\tau_o}(dt\wedgedot\etavector),$ where $c$ is a constant
  satisfying
  \begin{equation}
    \label{c}%
    \glossary{$c$}%
    c^2\int_{-\infty}^{\infty} e^{-\delta \tau_o}dt=2.
  \end{equation}
%   and $dt\wedgedot\etavector=((dt\wedge b)^+ , dt\, .\, \phi)$ for
%   $\etavector=(b,\phi),$ where in the second component dot denotes
%   Clifford multiplication. 
  (Recall that we identified spinors on $Y$
  with plus-spinors on the cylinder $\bbR \times Y.$ Clifford
  multiplication by $dt$ is just an isometry between plus- and
  minus-spinors on the cylinder.) Using this isomorphism, we can
%   find an orthonormal base $\overline{h_i}$ for
%   $\H^2_{-\delta}(\bbR\times Y)$ and
  express $$\Pi^2_o(\zetavector)=
  c^2\sum_{i=1}^{n}e^{-\delta\tau_o}(dt\wedgedot h_i)
  \int_{\bbR\times \displaystyle Y}\ip{\zetavector}{dt\wedgedot h_i}$$ for
  $\zetavector\in\E^2_{\delta}(\bbR\times Y).$ Using the facts above,
  a calculation shows that
% \note{say more}
  $$\tilde{\Pi}^2\circ\tilde{\Pi}^2\zeta=
  \frac{c^2}{2}\left(\int_{\bbR}\mu_o^2e^{-\delta\tau_o}\right)
  \tilde{\Pi}^2\zeta.$$
  Now, this formula, equation~\ref{c} and the fact that $\mu_o=1$
  on $C_{L_0}$ give the desired estimate.
  \\
  Finally, the right inverse can be constructed as follows. Choose a
  cut-off function $\beta,$ depending only on $t$ and supported in
  $C_1,$ such that $\int_{\bbR}\beta(t)dt=2.$ Then define the
  right-inverse $F$ by
  $F(\zetavector)=\frac{\beta}{c^2}e^{\delta\tau_o}\zetavector.$
  Since $\mu_o=1$ on the support of $\beta,$ it is easy
% \note{verify}
  to see that $\tilde{\Pi}^2(F\zetavector)=\zetavector$ for all
  $\zetavector.$ Clearly, $\|F\|$ is independent of $\ell$ and $L_0.$
\end{proof}
As a result, for $L_0$ large enough,
$\im(\tilde{\Pi}^2)\cap \ker(\tilde{\Pi}^2)=0$
and we obtain a decomposition
% \note{say more here}%
$\E^2_{-\delta}(X_\ell)=\im(\tilde{\Pi}^2) \plus \ker(\tilde{\Pi}^2).$
To see this, fix $L_0$ to satisfy $K^2_2 e^{-\delta L_0} <\frac12$ (or
any other number less than one, for that matter). Then, if $z\in
\im(\tilde{\Pi}^2)\cap \ker(\tilde{\Pi}^2),$ then one can express $z$ as
$z=\tilde{\Pi}^2(Fz).$ Thus,
$$
  \|z\|=\|z-\tilde{\Pi}^2 z\|=
  \|\tilde{\Pi}^2(Fz)-\tilde{\Pi}^2(\tilde{\Pi}^2(Fz))\|
  \leq K_2 e^{-\delta L_0}\|(Fz)\| \leq K_2^2 e^{-\delta L_0}\|z\|
  <\frac12\|z\|,
$$
which can not happen unless $z=0.$
\\
By the way, the above argument also shows that if $z\in
\im(\tilde{\Pi}^2),$ then $\|z-\tilde{\Pi}^2 z\|\leq\frac12\|z\|.$
Therefore, for such a $z,$ $\|z\|\leq 2\|\tilde{\Pi}^2 z\|.$ This will
be used below in the proof of lemma~\ref{Pi2}.
\\
Finally, the decomposition results from the fact that
$\im(\tilde{\Pi}^2)$ is finite-dimensional, being identified with
$\im(\Pi^2_o)=\H^2_{\delta}(\bbR\times Y).$
\par
Define a \textit{projection}
$\Pi^2:\E^2_{-\delta}(X_\ell) \cdr \E^2_{-\delta}(X_\ell)$ onto
$\im(\tilde{\Pi}^2)$ corresponding to this decomposition. Thus
$\im(\Pi^2) = \im(\tilde{\Pi}^2)$ and $\ker(\Pi^2) = \ker(\tilde{\Pi}^2).$
\begin{lemma}
  \label{Pi2}%
  \glossary{$K'_2$}%
  $\|\Pi^2-\tilde{\Pi}^2\|_{1,-\delta}\leq K'_2 e^{-\delta L_0}.$
%   \endpf
\end{lemma}
\begin{proof}
Decompose $\zetavector=z+z_0,$ where $z\in \im(\tilde{\Pi}^2)$ and
$z_0\in \ker(\tilde{\Pi}^2).$ Then
$$\Pi^2\zetavector-\tilde{\Pi}^2 \zetavector=z-\tilde{\Pi}^2 z=
\tilde{\Pi}^2(Fz)-\tilde{\Pi}^2(\tilde{\Pi}^2(Fz)).$$
Therefore, using lemma~\ref{tildePi2},
$$\|\Pi^2\zetavector-\tilde{\Pi}^2 \zetavector\|\leq
K_2 e^{-\delta L_0}\|Fz\|\leq K_2^2 e^{-\delta L_0}\|z\|\leq
K_2^2 K_3 e^{-\delta L_0}\|\zetavector\|,$$
where in the last inequality we have used the following remark
(\ref{K3}).
\end{proof}
\begin{rem}
\label{K3}%
\glossary{$K_3$}%
  If $\zetavector=z+z_0$ is a decomposition of $\zetavector,$ where
  $z\in \im(\tilde{\Pi}^2)$ and $z_0\in \ker(\tilde{\Pi}^2),$ then there
  is a constant $K_3$ such that $\|z\|+\|z_0\|\leq
  K_3\|\zetavector\|.$
\end{rem}
\noindent\textit{Proof.}
  This is a subsequence of the fact
  alluded to earlier.  Namely, we have
  $$\|z\|\leq 2\|\tilde{\Pi}^2 z\|= 2\|\tilde{\Pi}^2\zetavector\|\leq
  2\|\tilde{\Pi}^2\|.\|\zetavector\|,$$
  $$\|z_0\|\leq\|\zetavector\|+\|z\|\leq
  (1+2\|\tilde{\Pi}^2\|)\|\zetavector\|. \endpf $$
\par
We will explicitly identify $\im(\tilde{\Pi}^2)=\im(\Pi^2)$ with the
Zariski tangent space of $\M^3(Y)$ (at $[\etapt]),$ which is a
finite-dimensional vector space.
% \note{$\im(\Pi^2)$?}
% Recall that $\im(\Pi^0_2)=H^2(\E({\xipt}^0))$
% can be identified with $H^2(\E_Y(\etapt)).$
\begin{lemma}
  \label{identification}%
  \glossary{$\imath$}\glossary{$K_{\imath}$}%
% \note{$\imath$ defined on p. 64 of \cite{mm}}
  The linear map
  $$ \imath: H^1(\E(Y)|_\etapt) \cdr \im(\Pi^2) \subset
  \E^2_{-\delta}(X_\ell)$$ given by
  $$ \imath(\etavector)=c \mu_o
  e^{-\delta(\tau_o-\ell/2)}(dt\wedgedot\etavector)$$ is
  an isomorphism. It approaches an isometry as $L_0 \to \infty.$ More
  precisely, we have the following estimate. For some constant $K_{\imath},$
  $$(1-K_{\imath}e^{-\delta L_0/2})\leq \|\imath\|
  \leq(1+K_{\imath}e^{-\delta L_0/2}).$$   % \endpf$$
% [Is $\ell$ correct or $L_0?$]
%  \endpf $$
\end{lemma}
\par
We will later re-scale $\imath$ to fit it into an ``almost-commutative''
% quasi-commutative [?]
diagram. The last statement in the lemma will be used for estimating
$\|\imath\|$ in $\|\R_2\|$ (see proposition~\ref{R2}).
\begin{proof}
  This is a straightforward estimate. Only note that on
  $C_\ell\subset X_\ell,$ we are using $\tau_\ell,$ while
  $\tau_o=\ell-\tau_\ell$ is used on $X^o=\bbR\times Y.$ Comparing the
  two norms, therefore, we have
  $\|\etavector\|_{-\delta}=e^{-\delta\ell/2}\|\etavector\|_{\delta},$
  where the first norm is measured on $X_\ell$ and the second on
  $X^o.$ We will also use the fact that $\mu_o=1$ on $C_{L_0}.$
\end{proof}
\par % \vspace{1cm}
\bigskip
\glossary{\R_1}\glossary{\R_2}%
Now we head for constructing a right inverse for $\D^1.$ This will be
done in three steps. In the next subsection, we will construct a right
inverse $\R_1$ for $\D^1$ on the finite-codimensional subspace
$\ker(\Pi^2).$ In section~\ref{sec:R2}, we
construct a right inverse $\R_2$ for $\D^1$ on the transversal subspace
$\im(\Pi^2).$ There we use a stronger assumption that the fiber product of
the moduli spaces of the \cend manifolds is smooth in (a \nbhd of) the
point under consideration. This will be explained in more detail.
Finally, We show how to deform $\R_1+\R_2$ to get a right inverse $\R$
for $\D^1$ on all of $\E^2_{-\delta}(X_\ell)=\im(\Pi^2)\plus\ker(\Pi^2).$
% \note{Draw picture.}
%
\subsection{Right-Inverting $\D^1$ on Finite-Codimensional $\ker(\Pi^2)$}
% \label{sec:2.finite-codim-inverse}
\label{sec:R1}%
\glossary{$\R_1$}%
Recall that we glued the three operators $\R_+,$ $\R_-$ and $\R_o$ to
obtain $\tilde{\R}.$ We will be slightly modifying this operator to
establish the existence of a right inverse $\R_1$ for $\D^1$ on
$\ker(\Pi^2).$ We extend $\R_1$ by 0 on the complementary subspace
$\im(\Pi^2).$
\begin{prop}
  \label{R1}%
  \glossary{$R_1$}\glossary{$C_1$}\glossary{$N_1$}%
%   Suppose that
%   $([\xipt^+],[\xipt^-])\in\M(X^+,P^+)\times_U\M(X^-,P^-)$ and that
%   $[\xipt^\pm]$ are smooth points of their moduli spaces. Then if
  If $L_0$ is chosen sufficiently large, there is a constant $C_1$
\glossary{$C_1$}%
  such that the following holds. For all $\ell>4L_0,$
% Then there are constants $C_6,L_6>0$ such that the following holds for all
% $L_0\ge L_3$ and for all $L\ge 4L_0$
  there is an operator
  $$\R_1=\R_1(\xipt^+,\xipt^-,\ell):\E^2_{-\delta}(X_\ell) \to
  \E^1_{-\delta}(X_\ell),$$ such that
\begin{enumerate}
\item For all $\zetavector\in
  \ker(\Pi^2)\subset\E^2_{-\delta}(X_\ell),$
  $$({\mathrm I}-\Pi^2)\D^1\R_1(\zetavector)=\zetavector.$$
\item For all $\zetavector \in \im(\Pi^2)$ we have
  $$\R_1(\zetavector)=0.$$
\item The operator norm of $\R_1$ is bounded by $C_1,$ independent of
  $\ell$ and $L_0.$
\item Define $N_1=N_1(\xipt^+,\xipt^-,\ell)$ by setting
  $N_1=\Pi^2\D^1\R_1.$ Then $N_1^2=0$ and the norm of this operator
  satisfies $$\|N_1\|\leq\frac{C_1}{L_0}.$$
\item Moreover, $\R_1$ is asymptotically close to $\tilde{\R}:$
  $$\|\R_1-\tilde{\R}\|\leq\frac{C_1}{L_0}.$$
\end{enumerate}
% \endpf
% \note{[MM], pp. 72-74.}
\end{prop}
\begin{proof}
  \glossary{$J_1$}\glossary{$j_1$}%
  We are going to estimate the norm of the following operator, restricted to
  $\ker(\Pi^2),$
  $$({\mathrm I}-\Pi^2)\D^1\tilde{\R}:\ker(\Pi^2)\to \ker(\Pi^2).$$
%   , to show that it can be made close to identity. Therefore it would
%   have an inverse $J$ and then we would define
%   $\R_1=\tilde{\R}\circ J.$
%   \\
  Let $\zetavector\in \ker(\Pi^2)$ and decompose
  $$(\D^1\tilde{\R}-{\mathrm I})(\zetavector)=
    ({\mathrm I}-\Pi^2)(\D^1\tilde{\R}\zetavector-\zetavector)+
    \Pi^2(\D^1\tilde{\R}\zetavector),$$
  where the terms on the right hand side belong to $\ker(\Pi^2)$ and
  $\im(\Pi^2),$ respectively. Thus, according to remark~\ref{K3}, we
  have
  $$\|({\mathrm I}-\Pi^2)(\D^1\tilde{\R}\zetavector-\zetavector)\|\leq
    K_3 \|(\D^1\tilde{\R}-{\mathrm I})(\zetavector)\|.$$
%   Now, the element on the left side of this inequality is
  Now, in this norm comparison inequality, the element on the left
  side is 
  $({\mathrm I}-\Pi^2)\D^1\tilde{\R}\zetavector-\zetavector$ and the one
  on the right is
  $(\D^1\tilde{\R}-{\mathrm I}+\tilde{\Pi}^2)(\zetavector),$ whose norm,
  by lemma~\ref{first-estimates}, is bounded by
  $\frac{\tilde{K}}{L_0}\|\zetavector\|.$ Thus, we obtain
  $\|({\mathrm I}-\Pi^2)\D^1\tilde{\R}\zetavector-\zetavector\|\leq
  \frac{K_3\tilde{K}}{L_0}\|\zetavector\|$ and
  $$\|({\mathrm I}-\Pi^2)\D^1\tilde{\R}-{\mathrm I}\|\leq
  \frac{K_3\tilde{K}}{L_0}.$$
  Choose $L_0$ such that $\frac{K_3\tilde{K}}{L_0}<\frac12.$ Then the
  operator introduced at the beginning of the proof has an inverse of
  the form $J_1={\mathrm I}+j_1,$ where $j_1:\ker(\Pi^2)\to \ker(\Pi^2)$ satisfies
  $$\|j_1\|\leq\frac{2K_3\tilde{K}}{L_0}.$$
  Now set $\R_1=\tilde{\R}\circ J_1$ and extend $\R_1$ by zero on
  $\im(\Pi^2).$ The first three items are now immediate. To prove the
  fourth, notice that we only need to work on $\ker(\Pi^2),$ since
  $\R_1$ vanishes on $\im(\Pi^2).$ For $\zetavector\in \ker(\Pi^2),$ we
  have
  $$N_1\zetavector=\Pi^2\D^1\R_1(\zetavector)=
  \D^1\R_1(\zetavector)-\zetavector=
  \D^1\tilde{\R}J_1\zetavector-J_1\zetavector+j_1\zetavector.$$
  The last term already satisfies the desired estimate. On
  $\ker(\Pi^2),$ we also have
  $$
  \begin{array}[t]{rcl}
    \|\D^1\tilde{\R}J_1-J_1\|\leq
    \|\D^1\tilde{\R}-{\mathrm I}+\tilde{\Pi}^2\|.\|J_1\| & \leq &
    \frac{\tilde{K}}{L_0}(1+\frac{2K_3\tilde{K}}{L_0}) \\
      & \leq & \frac{2\tilde{K}}{L_0},
  \end{array}
  $$
  which establishes the fourth estimate. \\
  Finally, $\R_1-\tilde{\R}=\tilde{\R}J_1-\tilde{\R}=\tilde{\R}j_1$
  has the desired decay.
\end{proof}
\subsection{Right-Inverting $\D^1$ on Finite-Dimensional $\im(\Pi^2)$}
\label{sec:R2}%
\glossary{$(D)$}%
Consider the following diagram $(D).$
% \smallskip
$$
\label{D}%
\xymatrix{
H^1(\E_{-\delta}(X^+)|_{\xipt^+})\plus H^1(\E_{-\delta}(X^-)|_{\xipt^-})
\ar[r]^(0.70){r} \ar[d]_{\jmath} &
H^1(\E(Y)|_{\etapt}) \ar@{.>}@<1ex>[l]^(0.30){\rho} \ar[d]^{\imath}
\\
\E^1_{-\delta}(X_\ell) \ar[r]^{\D^1} & \E^2_{-\delta}(X_\ell).
}
\put(60,-25){$(D)$}%
%\eqno{(D)}\label{D}
$$
In this diagram,
$$
\glossary{$r$}%
r(\xivector_+,\xivector_-)=r_+(\xivector_+)-r_-(\xivector_-),
$$
$$
\glossary{$\jmath$}\glossary{$\jmath$ in [S] = $\imath_1$ in [MM]}%
\jmath(\xivector_+,\xivector_-)=\nu_+\xivector_+ +\nu_-\xivector_-,
$$
where $\nu_+$ and $\nu_-$ are certain cut-off functions, to be defined
below. First, define
$$
\glossary{$\nu$}%
\nu=\int_{-2L_0}^{2L_0}e^{-\delta \tau_o(s)}\mu_o(s)ds.
$$
Then, for $(t,y)\in\bbR\times Y,$
$$
\glossary{$\nu_-$}%
\nu_-(t,y)=\frac{1}{\nu}\int_{-2L_0}^{t}e^{-\delta\tau_o(s)}\mu_o(s)ds,
$$
$$
\glossary{$\nu_+$}%
\nu_+(t,y)=\frac{1}{\nu}\int_{t}^{2L_0}e^{-\delta\tau_o(s)}\mu_o(s)ds.
$$
We have $\nu_++\nu_-=1$ and these cut-off functions are constant
outside $C_{2L_0}.$ We will show shortly (in lemma~\ref{j}) that $\jmath$
is a quasi-isometry.
% \note{correct?}%
% close to an isometry onto its image for large $L_0.$
Also recall that the embedding $\imath,$ whose image identified with
$\im(\Pi^2),$ was defined in lemma~\ref{identification} by
$$
\glossary{$\imath$}%
\imath(\etavector)=
c\mu_o e^{-\delta(\tau_o-\ell/2)}(dt\wedgedot\etavector)
$$
and $\D^1:\E^1_{-\delta}(X_\ell) \to \E^2_{-\delta}(X_\ell)$ is the
differential of $\SW$ at
$\tilde{\xipt}_\ell=(\tilde{A}_\ell, \tilde{\Psi}_\ell)$ given by the
matrix
$$
\D^1=
\left(
  \begin{array}{cc}
    d^+ & -Dq|_{\tilde{\Psi}_\ell} \\
    .\frac12 \tilde{\Psi}_\ell  & \dirac_{\tilde{A}_\ell}
  \end{array}
\right) .
$$
\par
$\rho$ is a right inverse for $r$ and is essential in our
construction. Of course, to ensure the existence of such a right
inverse, we need the following 
\par
% \smallskip
\textbf{Transversality Assumption.}
\index{transversality assumption}%
The limiting maps
$\partial_+:\M^*(X^+)\rightarrow\M^{irr}(Y)$ and
$\partial_-:\M^*(X^-)\rightarrow\M^{irr}(Y)$ are
transversal at $\xipt^+$ and $\xipt^-,$ where
$\partial_+(\xipt^+)=\partial_-(\xipt^-)=\etapt.$ 
% =\xipt_\infty \\
In other words, the fiber product $\M^*(X^+)\times_U\M^*(X^-)$ is
smooth. Equivalently, the linear map
$$
\begin{array}{c}
r:H^1(\E_{-\delta}(X^+)|_{\xipt^+})\plus H^1(\E_{-\delta}(X^-)|_{\xipt^-})\to
H^1(\E(Y)|_{\etapt}) \\ r(\xi^+,\xi^-)=r_+(\xi^+)-r_-(\xi^-)
\end{array}
$$
is onto. $r_+$ and $r_-$ are the linearized versions of $\partial_+$ and
$\partial_-,$ respectively.
\par
Unfortunately, the diagram $(D)$ is not commutative; fortunately, it is
close to one, in the sense of the following lemma.
\begin{lemma}
\glossary{$K_D$}\glossary{$c_\ell$}%
  In diagram $(D),$ if $L_0$ is chosen large enough, there is a
  constant $K_D$ such that for all $\ell\geq 4L_0$ we have
% \note{\cite{mm}, p. 80}
  $$\|(\D^1\circ\jmath)-c_\ell.(\imath\circ r)\|_{1,-\delta}\leq
  K_D e^{-(\frac{\kappa+\delta}{2})(\ell-2L_0)},$$
  where $c_\ell=-\frac{1}{c\nu}e^{-\delta\ell/2}$ is a re-scaling
  factor.     % \endpf
\end{lemma}
\begin{proof}
% \note{cf. \cite{mm}, pp. 77, 61.}
% \note{$\jmath$ in [S] $=\imath_1$ in \cite{mm}}
\glossary{$\xivector_\pm=(a_\pm,\psi_\pm)$}%
  We start by computing $\D^1\circ\jmath(\xivector_+,\xivector_-),$
  where $\xivector_\pm=(a_\pm,\psi_\pm).$ Note that this is supported
  on $C_{2L_0}.$ A component-wise calculation shows
  $$\D^1(\nu_+\xivector_+)=d\nu_+\wedgedot\xivector_+ +\nu_+\D^1\xivector_+$$
  and the same equation holds for $\D^1(\nu_-\xivector_-)$ with all
  the $+$'s and $-$'s reversed. Now,
  $$
  \begin{array}{rcl}
  \D^1(\nu_+\xivector_+)
  &=&
  \frac{-1}{\nu}e^{-\delta\tau_o}\mu_o dt\wedgedot\xivector_+
  +\nu_+\D^1\xivector_+
  \\
  &=&
  \frac{-1}{\nu}e^{-\delta\tau_o}\mu_o dt\wedgedot r_+(\xivector_+)
  +\frac{-1}{\nu}e^{-\delta\tau_o}\mu_o
  dt\wedgedot(\xivector_+-r_+(\xivector_+))
  +\nu_+\D^1\xivector_+.
  \end{array}
  $$
  We also get a similar formula for $\D^1(\nu_-\xivector_-),$ except
  for the sign of the first two terms. Therefore,
  $$
  \begin{array}[t]{rcl}
  (\D^1\circ\jmath)(\xivector_+,\xivector_-)
  &=&\frac{-1}{\nu}e^{-\delta\tau_o}\mu_o
     dt\wedgedot r(\xivector_+,\xivector_-) \\
  & &+\frac{-1}{\nu}e^{-\delta\tau_o}\mu_o
     \Bigl(dt\wedgedot(\xivector_+-r_+(\xivector_+))
          -dt\wedgedot(\xivector_+-r_+(\xivector_+))\Bigr)\\
  & &+\nu_+\D^1\xivector_+ +\nu_-\D^1\xivector_-, % \\
  \end{array}
  $$
  in which the first term on the right hand side is just
  $c_\ell.(\imath\circ r)(\xivector_+,\xivector_-).$ Thus we
  have obtained
  $$
  \begin{array}{rl}
  \Bigl((\D^1\circ\jmath)-c_\ell.(\imath\circ r)\Bigr)
        (\xivector_+,\xivector_-)
  =& \frac{-1}{\nu}e^{-\delta\tau_o}\mu_o
     \Bigl(dt\wedgedot(\xivector_+-r_+(\xivector_+))
          -dt\wedgedot(\xivector_+-r_+(\xivector_+))\Bigr) \\
   & +\nu_+\D^1\xivector_+ + \nu_-\D^1\xivector_-
  \end{array}
  $$
  Now the result follows from the exponential decay of solutions on
  $X_\ell^\pm$
  $$\|\xivector_\pm-r_+(\xivector_\pm)\|\leq
  e^{-\frac{\kappa}{2}(\ell-2L_0)}\|\xivector_\pm\|$$
  and the fact that $\D^1\xivector_\pm$ can be expressed in terms of
  the partition of unity $\lambda,$ the components of $\xipt_\pm,$
  $\xivector_\pm,$ and their derivatives.
\end{proof}
\par
Here is the analog of lemma~\ref{identification}.
\begin{lemma}
% \note{\cite{mm}, pp. 77, 61.}
% \note{$\jmath$ in [S] $=\imath_1$ in \cite{mm}}
\label{j}%
\glossary{$\jmath$}\glossary{$K_\jmath$}%
  The linear map
  $$\jmath:H^1(\E_{-\delta}(X^+)|_{\xipt^+})\plus
  H^1(\E_{-\delta}(X^-)|_{\xipt^-}) \to \E^1_{-\delta}(X_\ell)$$
  in diagram (D) defined by
  $\jmath(\xivector_+,\xivector_-)=\nu_+ \xivector_+ + \nu_- \xivector_-$
  is a quasi-isometry, satisfies
  $$(1-K_{\jmath}e^{-\delta L_0})\leq \|\jmath\| \leq
  (1+K_{\jmath}e^{-\delta l_0})$$
  and approaches an isometry as $L_0 \to \infty.$ Here $\xivector_+$
  and $\xivector_-$ are harmonic representatives of the corresponding
  cohomologies.       % \endpf
\end{lemma}
\begin{proof}
\glossary{$\jmath_1$}%
  Let's first define $$\jmath_1:H^1(\E_{-\delta}(X^+)|_{\xipt^+})\plus
  H^1(\E_{-\delta}(X^-)|_{\xipt^-}) \to \E^1_{-\delta}(X_\ell)$$
  by
  $\jmath_1(\xivector_+,\xivector_-)=\mu_+\xivector_++\mu_-\xivector_-.$
  It is straightforward to see that $\jmath_1$ is an isomorphism onto
  $\im(\tilde{\Pi}^1)$ and satisfies
  $1-K_{\jmath_1}e^{-\delta L_0/2}\leq \|\jmath_1\|\leq
  1+K_{\jmath_1}e^{-\delta L_0/2}.$ Moreover, a calculation shows
  that there is a constant $K_{\jmath}$ such that for all $L_0 \gg 0$
  and all $\ell> 4L_0,$ we have $\|\jmath-\jmath_1\|\leq
  K_{\jmath} e^{-\delta L_0}.$ Thus we get the desired estimate
  on $\|\jmath\|.$
\end{proof}
We are now in a position to state the main result of this section,
which is the counterpart of proposition~\ref{R1}.
\begin{prop}
\label{R2}%
\glossary{$C_2$}\glossary{$R_2$}%
%   Suppose that $([\xipt^+],[\xipt^-])$ is a smooth point of the fiber
%   product ${\cal M}(X^+,P^+)\times_U{\cal M}(X^-,P^-).$
  Suppose that the \emph{transversality assumption} holds for
  $\xipt^+$ and $\xipt^-.$
  Then, if $L_0$ is chosen large enough, there is a constant $C_2$ such
  that the following holds. For all $\ell>4L_0,$ there is an operator
  $$\R_2=\R_2(\xipt^+,\xipt^-,\ell): \E^2_{-\delta}(X_\ell) \to
  \E^1_{-\delta}(X_\ell),$$ 
% \im(\imath)\to \im(\jmath),$$ %$$
  such that
  \begin{enumerate}
  \item For all $\zeta \in \im(\Pi^2) % =\im(\imath)
    \subset \E^2_{-\delta}(X_\ell),$
    $$\Pi^2 \D^1 \R_2(\zeta)=\zeta.$$
  \item For all $\zeta \in \ker(\Pi^2)$ we have $$\R_2(\zeta)=0.$$
  \item The operator norm of $\R_2$ satisfies $$\| \R_2 \| \leq C_2
    e^{\delta\ell/2}.$$
  \item Define $N_2=N_2(\xipt^+,\xipt^-,\ell)$ by setting
    $N_2=({\mathrm I}-\Pi^2)\D^1\R_2.$ Then $N_2^2=0$ and the norm of this
    operator satisfies
    $$\|N_2\|\leq C_2e^{-(\kappa-\delta)L_0}.$$
  \end{enumerate}
\end{prop}
\begin{proof}
\glossary{$\tilde{\R}_2$}\glossary{$\tilde{C}_2$}\glossary{K_D}\glossary{$J_2$}%
  First, using the almost-commutative diagram $(D),$ define
  $$\tilde{\R}_2:=\frac{1}{c_\ell}(\jmath\circ\rho\circ\imath^{-1}):
  \im(\imath)=\im(\Pi^2)\to \im(\jmath).$$
  Then, it is easy to see that
  $\|\tilde{\R}_2\|\leq\tilde{C}_2 e^{\delta\ell/2}$ for some
  constant $\tilde{C}_2$ and that
  $$\|\Pi^2\circ\D^1\circ\tilde{\R}_2-{\mathrm I}\|\leq K_D
  e^{-(\frac{\kappa+\delta}{2})(\ell-2L_0)}e^{\delta\ell/2}\|\rho\|
  .\|\imath^{-1}\|.$$
  Thus one can choose $L_0$ sufficiently large so that for all
  $\ell>4L_0,$ we have
  $$\|\Pi^2\circ\D^1\circ\tilde{\R}_2-{\mathrm I}\|\leq\frac12$$
  and therefore $\Pi^2 \circ \D^1 \circ \tilde{\R}_2$ has an inverse
  $J_2$ of norm at most 2. Setting $\R_2=\tilde{\R}_2 \circ J_2$
  (and extending by zero to the complement of $\im(\Pi^2)$ in
  $\E^2_{-\delta}(X_\ell)$) gives the desired right inverse.
\end{proof}
\par
Notice that the proposition above implies that
$\D^1\circ\R_2={\mathrm I}+N_2$ on $\im(\Pi^2),$ where $N_2$ is a
nilpotent operator.  Proposition~\ref{R1} implied a similar statement,
that $\D^1\circ\R_1={\mathrm I}+N_1$ on $\ker(\Pi^2),$ with a nilpotent
$N_1.$ Both $N_1$ and $N_2$ depend on $\ell,$ as well as on $\xipt^+$
and $\xipt^-,$ since $\R_1$ and $\R_2$ do.
\subsection{Gluing Monopoles}
\label{sec:monopolesglue}%
\glossary{$\R_0$}%
It is easy now to construct a full right inverse for $\D^1.$ To
begin with, set $\R_0:=\R_1+\R_2.$ For a sufficiently large $L_0,$
both of the propositions~\ref{R1} and \ref{R2} hold and we have
$$\D^1 \circ \R_0 ={\mathrm I}+N_1+N_2.$$
(To check this identity, verify it on elements of $\ker(\Pi^2)$ and
$\im(\Pi^2).$) \\
Since $N_1^2=N_2^2=0,$
$$({\mathrm I}+N_1+N_2)({\mathrm I}-N_1-N_2)={\mathrm I}-N_1 N_2-N_2 N_1$$
and ${\mathrm I}+N_1+N_2$ has an inverse $J$ whenever
$2\|N_1\|.\|N_2\|<1.$ Moreover,
$$\|J\| \leq \frac{1+\|N_1\|+\|N_2\|}{1-2\|N_1\|.\|N_2\|}.$$
\glossary{$J$}%
But $\|N_1\|.\|N_2\|$ is bounded by $\frac{C_1C_2}{L_0}
e^{-(\kappa-\delta)L_0},$ which can be made arbitrarily small by
choosing $L_0$ large enough, since we chose $\delta < \kappa.$ So an
inverse $J$ with $\|J\|<2$ exists and we define $\R:=\R_0 \circ J.$
Then $\R$ is a right inverse for $\D^1$ and we have
\begin{prop}
\label{R}%
\glossary{$R$}\glossary{$C$}%
%   Suppose that $([\xipt^+],[\xipt^-])$ is a smooth point of the fiber
%   product ${\cal M}(X^+,P^+)\times_U{\cal M}(X^-,P^-).$
  Suppose that the \emph{transversality assumption} holds for
  $\xipt^+$ and $\xipt^-.$
  Then the operator
  $\D^1=D(\SW):\E^1_{-\delta}(X_\ell)\to\E^2_{-\delta}(X_\ell)$
  has a right inverse
  \index{right inverse}
  $$\R=\R_\ell(\xipt^+,\xipt^-):
  \E^2_{-\delta}(X_\ell)\to\E^1_{-\delta}(X_\ell),$$
  for each $\ell>4L_0,$ $L_0\gg 1$, whose norm satisfies
  $\|\R\|\leq Ce^{\delta\ell/2}.$ \endpf
\end{prop}
\par
{\bf Remark.} A review of the statements of this section shows that if
$U\subset\M(Y)$ is an open set whose closure contains no reducible
points, then, in each statement, we can replace $\kappa$ by an
exponent $\kappa(U)$ which works for all $\etapt\in U.$ This includes,
in particular, propositions \ref{exp-decay}, \ref{R1}, \ref{R2} and
\ref{R}. Therefore, it makes sense to consider the derivatives of the
operators in question and estimate their norms. It is not hard to see
that, in each case, the derivatives decay (or grow) exponentially with
the same exponent as the operators themselves.
\par
It is time to introduce our contraction mapping, whose fixed point
would be the correcting perturbation term for our approximately-glued
monopole.
\begin{lemma}
\label{fixed-point}%
\glossary{$\zetahat=\zetahat_\ell(\xipt^+,\xipt^-)$}%
\glossary{$\F=\F_\ell(\xipt^+,\xipt^-)$}\glossary{$C'$}%
  The following self-map $\F=\F_\ell(\xipt^+,\xipt^-)$ of the Hilbert
  space $\E^2_{-\delta}(X_\ell)$
% $=\Omega^2_{+,1,-\delta}(X_\ell;i\bbR)\plus\Gamma_{1,-\delta}(S^-)$
  is a contraction
\index{contraction map}%
  on a ball for $\ell>4L_0,$ $L_0\gg 1,$ and therefore has a unique
  fixed point
%  \index{fixed point}
  $\zetahat=\zetahat_\ell(\xipt^+,\xipt^-).$
% =(\hat\alpha,\hat s).
  $$\F:\E^2_{-\delta}(X_\ell)\to\E^2_{-\delta}(X_\ell)=
  \Omega^2_{+,1,-\delta}(X_\ell;i\bbR)\plus\Gamma_{1,-\delta}(S^-)$$
  $$\F(\zetavector)=  % =\F(\alpha,s)=
  -\SW(\tilde{\xipt}_\ell)+(q(\psi),-\frac12 a.\psi),$$
  where $(a,\psi)=\R(\zetavector)$ and % =\R(\alpha,s)$ and %$
  $$\R:\Omega^2_{+,1,-\delta}(X_\ell;i\bbR)\plus\Gamma_{1,-\delta}(S^-)
  \to\wforms12{-\delta}(X_\ell;i\bbR)\plus\Gamma_{2,-\delta}(S^+)$$
  is a right inverse of $\D^1$ as constructed before. \\
  Moreover, we have the following estimates on the norm of the fixed
  point and its image.
  $$\|\zetahat\|_{1,-\delta}\leq C'e^{-\kappa(\ell-2L_0)},$$
  $$\|\R(\zetahat)\|_{2,-\delta}\leq
  C'e^{-\kappa(\ell-2L_0)}e^{\delta L_0}.$$
  Furthermore, $\zetahat$ varies smoothly with $\xipt^+$ and
  $\xipt^-,$ so that if $\xipt^\pm(t)$ are smooth, one-parameter
  families in $\C(X^\pm)$ with the same irreducible limiting value,
  then $\zetahat(t)=\zetahat_\ell(\xipt^+(t),\xipt^-(t))$ is also a
  smooth one-parameter family and if $\zetahat'$ denotes
  $\frac{d}{dt}\zetahat(t)|_{t=0},$ then we have
  $$\|\R(\zetahat)'\|_{2,-\delta}\leq
  C'e^{-\kappa(\ell-2L_0)}e^{\delta L_0}
  \Bigl(\|(\xipt^+)'\|+\|(\xipt^-)'\|\Bigr).$$
\end{lemma}
{\bf Note.} The fact that the norm of $\R(\zetahat)'$
% $\|\R(\zetahat)'\|_{2,-\delta}$
is exponentially decreasing despite the possible exponential growth of
the operator norm of $\R$ is due to the quadratic nature of
$\F(\zetavector)$ in $\zetavector.$ This can be seen during the proof.
\glossary{G}\glossary{\tilde{c}_1}\glossary{\tilde{c}_2}%
\begin{proof}
  Let $B(0,R)$ denote the ball of radius $R$ around the origin in the
  Hilbert Space
  $\E^2_{-\delta}(X_\ell)=\Omega^2_{+,1,-\delta}(X_\ell;i\bbR) \plus
  \Gamma_{1,-\delta}(S^-).$ We are going to show that there is a
  constant $G$ such that for all $\ell\geq 4L_0,$ the restriction
  $\F|:B(0,Ge^{-2\delta\ell})\to B(0,Ge^{-2\delta\ell})$ 
  to the ball of radius $Ge^{-2\delta\ell}$ 
  is a $\frac12$-contraction.  First, we consider the norm of
  $\F(0)=-\SW(\tilde{\xipt}_\ell).$ By proposition~\ref{exp-decay},
  there is a constant $\tilde{C}$ such that
  $$\|\F(0)\|_{1,-\delta}=\|\SW(\tilde{\xipt}_\ell)\|_{1,-\delta}
  \leq \tilde{C}e^{-(\kappa+\frac{\delta}{2})(\ell-2L_0)}.$$
  Now,
  let's estimate each of the components of
\begin{equation}
\label{eq:F}%
\F(\zetavector_1)-\F(\zetavector_2)=
(q(\psi_1)-q(\psi_2),-\frac12(a_1.\psi_1-a_2.\psi_2)).
\end{equation}
  For the first component, we have
  $$q(\psi_1)-q(\psi_2)=\flat(\psi_1+\psi_2,\psi_1-\psi_2),$$
  where $\flat$ denotes the symmetric bilinear form associated to $q$
  and its point-wise norm is bounded above by
  $|\flat(\psi_1+\psi_2,\psi_1-\psi_2)|\leq
  2|\psi_1+\psi_2||\psi_1-\psi_2|.$ Therefore,
%   using the Cauchy-Schwarz ineqality
%   $$
% \note{see below}
\begin{eqnarray}%{rcl}
\label{eq:c1}%
  \|q(\psi_1)-q(\psi_2)\|_{1,-\delta} & \leq &
  2\, \|\psi_1+\psi_2\|_{L^4_{1,-\delta}}
  \|\psi_1-\psi_2\|_{L^4_{1,-\delta}} \nonumber\\
  & \leq &
  2\, e^{2\frac{\delta}{4}\ell}\|\psi_1+\psi_2\|_{L^4_{1,-2\delta}}
  \|\psi_1-\psi_2\|_{L^4_{1,-2\delta}} \nonumber\\ %??
  & \leq &
  \tilde{c}_1 e^{\delta\ell}\|\psi_1+\psi_2\|_{2,-\delta}
  \|\psi_1-\psi_2\|_{2,-\delta}.
%   \tilde{c}_1 e^{\delta\ell}\|\psi_1+\psi_2\|_{L^2_{2,-\delta}}
%   \|\psi_1-\psi_2\|_{L^2_{2,-\delta}}.
%   \begin{picture}(0,0)\put(52,0){(1)}\end{picture}
\end{eqnarray}
%   $$
  For the second component, we similarly write
%   $$
\begin{eqnarray}%{rcl}
\label{eq:c2}%
  \|a_1.\psi_1-a_2.\psi_2\|_{1,-\delta} & = &
  \|a_1(\psi_1-\psi_2)+(a_1-a_2)\psi_2\|_{1,-\delta} \nonumber\\
  & \leq &
  \|a_1\|_{L^4_{1,-\delta}}\|\psi_1-\psi_2\|_{L^4_{1,-\delta}}+
  \|\psi_2\|_{L^4_{1,-\delta}}\|a_1-a_2\|_{L^4_{1,-\delta}} \nonumber\\
  & \leq &
  e^{2\frac{\delta}{4}\ell}\bigl(
  \|a_1\|_{L^4_{1,-2\delta}}\|\psi_1-\psi_2\|_{L^4_{1,-2\delta}}+
  \|\psi_2\|_{L^4_{1,-2\delta}}\|a_1-a_2\|_{L^4_{1,-2\delta}}\bigr)
  \nonumber\\
  & \leq &
  \tilde{c}_2 e^{\delta\ell}\bigl(
  \|a_1\|_{2,-\delta}\|\psi_1-\psi_2\|_{2,-\delta}+
  \|\psi_2\|_{2,-\delta}\|a_1-a_2\|_{2,-\delta}\bigr).
%   \tilde{c}_2 e^{\delta\ell}\bigl(
%   \|a_1\|_{L^2_{2,-\delta}}\|\psi_1-\psi_2\|_{L^2_{2,-\delta}}+
%   \|\psi_2\|_{L^2_{2,-\delta}}\|a_1-a_2\|_{L^2_{2,-\delta}}\bigr).
\end{eqnarray}
%   $$
%   \|\|_{L^4_{1,-\delta}}
%   $$a_1.\psi_1-a_2.\psi_2=a_1(\psi_1-\psi_2)+(a_1-a_2)\psi_2.$$
  So, to finish the estimates on the components, we need estimates on
  the $L^2_{2,-\delta}$-norms of $a_i$ and $\psi_i,$ for $i=1,2,$ as
  well as the differences $a_1-a_2$ and $\psi_1-\psi_2.$ First, note
  that $(a_i,\psi_i)=\R(\zetavector_i)$ for $i=1,2,$ so that each of
  $\|a_i\|_{2,-\delta}$ and $\|\psi_i\|_{2,-\delta}$ is bounded above
  by
\begin{equation}
\label{eq:a-and-psi-estimate}%
  \|\R(\zetavector_i)\|_{2,-\delta}\leq
  Ce^{\delta\ell/2}\|\zetavector_i\|_{1,-\delta}.
\end{equation}
%   $$\|\R(\zetavector_i)\|_{2,-\delta}\leq
%   Ce^{\delta\ell/2}\|\zetavector_i\|_{1,-\delta}.$$
  Similarly, $a_1-a_2$ and $\psi_1-\psi_2$ are the two components of
  $\R(\zetavector_1)-\R(\zetavector_2),$ so that both of
  $\|a_1-a_2\|_{2,-\delta}$ and $\|\psi_1-\psi_2\|_{2,-\delta}$ are
  bounded above by
\begin{equation}
\label{eq:a-and-psi-difference-estimate}%
  \|\R(\zetavector_1)-\R(\zetavector_2)\|_{2,-\delta}\leq
  \|D_\zeta\R\|.\|\zetavector_1-\zetavector_2\|_{1,-\delta}\leq
  Ce^{\delta\ell/2}\|\zetavector_1-\zetavector_2\|_{1,-\delta}.
\end{equation}
%   $$\|\R(\zetavector_1)-\R(\zetavector_2)\|_{2,-\delta}\leq
%   \|\R'\|.\|\zetavector_1-\zetavector_2\|_{1,-\delta}\leq
%   Ce^{\delta\ell/2}.\|\zetavector_1-\zetavector_2\|_{1,-\delta}.$$
  These inequalities, in conjunction with estimates \ref{eq:c1} and
  \ref{eq:c2}, show that for $\zetavector_1,\zetavector_2$ in a ball
  $B(0,R)$ of radius $R,$ we have
  $$\|q(\psi_1)-q(\psi_2)\|_{1,-\delta}\leq 2\tilde{c}_1C^2R
  e^{2\delta\ell}\|\zetavector_1-\zetavector_2\|_{1,-\delta}$$
  $$\|a_1.\psi_1-a_2.\psi_2\|_{1,-\delta}\leq 2\tilde{c}_2C^2R
  e^{2\delta\ell}\|\zetavector_1-\zetavector_2\|_{1,-\delta}.$$
  Combining with equation \ref{eq:F}, we obtain
  $$\|\F(\zetavector_1)-\F(\zetavector_2)\|_{1,-\delta}\leq
  \tilde{c}R
%   (2\tilde{c}_1+\tilde{c}_2)C^2R
  e^{2\delta\ell}\|\zetavector_1-\zetavector_2\|_{1,-\delta}$$
  for some constant $\tilde{c}$ and $\F$ will be a
  $\frac12$-contraction for $R=\frac{1}{2\tilde{c}}e^{-2\delta\ell}.$
  \\
  Now, the unique fixed point of $\F$ can be obtained by finding the
  limit of the iterations of any point in the ball. Therefore, the
  sequence of iterations $\{\F^{\circ n}(0)\}$ converges to the fixed
  point $\hat\zetavector$ and we have
\begin{equation}
\label{eq:zetahat-estimate}%
  \|\hat\zetavector\|_{1,-\delta}\leq 2\|\F(0)\|_{1,-\delta}\leq
  2\tilde{C}e^{-(\kappa+\frac{\delta}{2})(\ell-2L_0)}.
\end{equation}
%   $$\|\hat\zetavector\|_{1,-\delta}\leq 2\|\F(0)\|_{1,-\delta}\leq
%   2\tilde{C}e^{-(\kappa+\frac{\delta}{2})(\ell-2L_0)}.$$
  From here, the estimates claimed in the theorem on
  $\|\hat\zetavector\|$ and $\|\R(\hat\zetavector)\|$ follow. \\
  Finally, to estimate the norm of
  $\R(\zetahat)'=\R'(\zetahat)+\R(\zetahat'),$ we need to estimate
  $\|\zetahat'\|_{1,-\delta}$ first. To avoid complicated notation, we
  will write $\zetahat$ for $\zetahat(t)|_{t=0}$ and $\zetahat'$ for
  $\frac{d}{dt}\zetahat(t)|_{t=0}.$ We will also consider the
  one-parameter family of operators $\F_t=\F(\xipt^+(t),\xipt^-(t))$
  and write $\F$ for $\F(0)$ and $\F'$ for the $t$-derivative of
  $\F_t$ at $t=0.$ Similar notation was already used in the case of
  $\R$ and $\R'$ at the beginning of this paragraph.\\
  By $t$-differentiating the fixed point equation
  $\F_t(\zetahat(t))=\zetahat(t)$ at $t=0$ we obtain
  $$D_\zeta\F(\zetahat')+\F'(\zetahat)=\zetahat'.$$
  (To see this, write $\F_t=\F+t\F'+o(t^2)$ and
  $\zetahat(t)=\zetahat+t\zetahat'+o(t^2),$ then substitute in the
  fixed point equation and find the coefficient of $t.$ Note that $\F$
  is not a linear map, so the coefficient of $t$ in
  $\F(t\zetahat')$ is $D_\zeta\F(\zetahat').$) \\
  We can rewrite the last equation as
  $({\mathrm I}-D_\zeta\F)(\zetahat')=\F'(\zetahat).$ Since $\F$ is a
  contraction, $\|D_\zeta\F\|\leq\frac12,$ so ${\mathrm I}-D_\zeta\F$ is
  invertible and $\|({\mathrm I}-D_\zeta\F)^{-1}\|\leq 2.$ Therefore,
\begin{equation}
\label{eq:zetahat'}%
  \|\zetahat'\|_{1,-\delta}\leq 2\|\F'(\zetahat)\|_{1,-\delta}.
\end{equation}
%   $$\|\zetahat'\|_{1,-\delta}\leq 2\|\F'(\zetahat)\|_{1,-\delta}.$$
  On the other hand, as in the preceding arguments,
% \note{verify carefully}
%   (cf. equations \ref{eq:F}, \ref{eq:c1}, \ref{eq:c2},
%   \ref{eq:a-and-psi-estimate} and
%   \ref{eq:a-and-psi-difference-estimate}),
  we find
  $$\|\F'(\zetahat)\|_{1,-\delta}\leq
  C'e^{\delta\ell}\|\R'(\zetahat)\|_{2,-\delta}
  \|\R(\zetahat)\|_{2,-\delta}
  +\|\SW(\tilde\xipt_\ell)'\|_{1,-\delta}.$$
  We moreover have
  $$\|\R'(\zetahat)\|_{2,-\delta}\|\R(\zetahat)\|_{2,-\delta}\leq
  C'e^{\delta\ell}\Bigl(\|(\xipt^+)'\|+\|(\xipt^-)'\|\Bigr)
  \|\zetahat\|_{1,-\delta}^2.$$
  Combining the last two inequalities and using the estimate on
  $\|\SW(\tilde\xipt_\ell)'\|_{1,-\delta}$
  (cf. proposition~\ref{exp-decay}) and the facts $\delta<\kappa$ and
  $\ell>4L_0,$  we see that, for some constant $C',$
  $$\|\F'(\zetahat)\|_{1,-\delta}\leq
  C'e^{-(\kappa+\delta/2)(\ell-2L_0)}
  \Bigl(\|(\xipt^+)'\|+\|(\xipt^-)'\|\Bigr)$$
  and, using (\ref{eq:zetahat'}), we get, for some $C',$
\begin{equation}
\label{eq:zetahat'-estimate}%
  \|\zetahat'\|_{1,-\delta}\leq C'e^{-(\kappa+\delta/2)(\ell-2L_0)}
  \Bigl(\|(\xipt^+)'\|+\|(\xipt^-)'\|\Bigr).
\end{equation}
  Now, using the equation $\R(\zetahat)'=\R'(\zetahat)+\R(\zetahat'),$
  the estimates on the norms of $\zetahat$ and $\zetahat'$ (equations
  \ref{eq:zetahat-estimate} and \ref{eq:zetahat'-estimate}) and the
  estimates on the operator norms $\|\R\|$ and $\|\R'\|$ (proposition
  \ref{R} and its following remark) gives the desired estimate on
  $\|\R(\zetahat)'\|_{2,-\delta}.$
\end{proof}
We are now in the final stage of gluing. Recall that we started with a
couple of solutions $\xipt^+$ and $\xipt^-$ on $X^+$ and $X^-,$
respectively. They solved $\SW(\xipt^\pm)=0.$
% $\SW_{h_+}(\xipt^+)=0$ and $\SW_{h_-}(\xipt^-)=0.$ 
Then we truncated and glued these solutions, using a partition of
unity, to obtain an \textit{approximate} solution
$\tilde{\xipt}_\ell.$ Set
$\xipt_\ell=\tilde{\xipt}_\ell+\R(\zetahat),$ where
$\R:\E^2(X_\ell)\to\E^1(X_\ell)$ is the right inverse of $\D^1$
constructed before (proposition~\ref{R}) and
$\zetahat=(\hat\alpha,\hat s)\in\E^2(X_{\ell})=
\Omega^2_{+,1}(X_\ell;i\bbR)\plus\Gamma_1(S^-)$ is the unique fixed
point of $\F$ of proposition~\ref{fixed-point}.  As we have
  $$\SW(\xipt+\xivector)=
  \SW(\xipt)+\D^1(\xivector)+(-q(\psi),\frac12a.\psi),$$
where $\xivector=(a,\psi),$ we obtain
  $\SW(\tilde{\xipt}_\ell+\R(\zetahat))=0,$
so $\xipt_\ell$ is a solution of $\SW.$
\endpf
% 
%
%%%%%%%%%%%%%%%%%%%%%%%%%%%%% End Text Here %%%%%%%%%%%%%%%%%%%%%%%%%%%
%
%
%
%%%%%%%%%%%%%%% Bibliography %%%%%%%%%%%%
%
\newpage
% \nocite{mm}
% \bibliography{bib,more}

\begin{thebibliography}{MCW1}
 
\bibitem[C]{c}
Weimin Chen.
\newblock The seiberg-witten theory of homology 3-spheres.
\newblock unpublished preprint,
\arxiv{dg-ga/9703009}$\!\!.$
 
\bibitem[FS]{fs}
Ronald Fintushel and Ronald~J. Stern.
\newblock Knots, links, and {$4$}-manifolds.
\newblock {\em Invent. Math.}, 134(2):363--400, 1998.
\MRhref{99j:57033}$\!\!.$
 
\bibitem[J]{joe}
Dosang Joe.
\newblock {\em Symplectic structures on connected sums with a ruled surface and
  product formulas for Seiberg-Witten invariants along a nilmanifold}.
\newblock PhD thesis, Columbia Univerity, 1998.
 
\bibitem[KM]{km}
P.~B. Kronheimer and T.~S. Mrowka.
\newblock The genus of embedded surfaces in the projective plane.
\newblock {\em Math. Res. Lett.}, 1(6):797--808, 1994.
\MRhref{96a:57073}$\!\!.$

\bibitem[MCW1]{mcw1}
Matilde Marcolli, Alan Carey, and Bai-Ling Wang.
\newblock Exact triangles in seiberg-witten floer theory. part i: the geometric
  triangle.
\newblock unpublished preprint,
\arxiv{math.DG/9907065}$\!\!.$
 
\bibitem[MW2]{mw2}
Matilde Marcolli and Bai-Ling Wang.
\newblock Exact triangles in seiberg-witten floer theory. part ii: geometric
  limits of flow lines.
\newblock unpublished preprint,
\arxiv{math.DG/9907080}$\!\!.$

\bibitem[M]{m}
John~W. Morgan.
\newblock {\em The {S}eiberg-{W}itten equations and applications to the
  topology of smooth four-manifolds}, volume~44 of {\em Mathematical
Notes}.
\newblock Princeton University Press, Princeton, NJ, 1996.
\MRhref{97d:57042}$\!\!.$
 
\bibitem[MMR]{mmr}
John~W. Morgan, Tomasz Mrowka, and Daniel Ruberman.
\newblock {\em The ${L}\sp 2$-moduli space and a vanishing theorem for
  {D}onaldson polynomial invariants}.
\newblock Monographs in Geometry and Topology, II. International Press,
  Cambridge, MA, 1994.
\MRhref{95h:57039}$\!\!.$
                                                                                
\bibitem[MM]{mm}
John~W. Morgan and Tomasz~S. Mrowka.
\newblock On the gluing theorem for instantons on manifolds containing long
  cylinders.
\newblock unpublished preprint.
 
\bibitem[MST]{mst}
John~W. Morgan, Zolt{\'a}n Szab{\'o}, and Clifford~Henry Taubes.
\newblock A product formula for the {S}eiberg-{W}itten invariants and the
  generalized {T}hom conjecture.
\newblock {\em J. Differential Geom.}, 44(4):706--788, 1996.
\MRhref{97m:57052}$\!\!.$
 
\bibitem[MOY]{moy}
Tomasz Mrowka, Peter Ozsv{\'a}th, and Baozhen Yu.
\newblock Seiberg-{W}itten monopoles on {S}eifert fibered spaces.
\newblock {\em Comm. Anal. Geom.}, 5(4):685--791, 1997.
\arxiv{math.GT/9612221}$\!\!.$
\MRhref{98m:58017}$\!\!.$
 
\bibitem[S]{s}
Pedram Safari.
\newblock {\em A gluing theorem for Seiberg-Witten moduli spaces}.
\newblock PhD thesis, Columbia Univerity, 2000.
 
\bibitem[T1]{tswgr1}
Clifford~Henry Taubes.
\newblock The {S}eiberg-{W}itten and {G}romov invariants.
\newblock {\em Math. Res. Lett.}, 2(2):221--238, 1995.
\MRhref{96a:57076}$\!\!.$
 
\bibitem[T]{tswgr}
Clifford~Henry Taubes.
\newblock {\em Seiberg {W}itten and {G}romov invariants for symplectic
  {$4$}-manifolds}, volume~2 of {\em First International Press Lecture
Series}.
\newblock International Press, Somerville, MA, 2000.
\newblock Edited by Richard Wentworth.
\MRhref{2002j:53115}$\!\!.$
 
\end{thebibliography}
\bibliographystyle{plain}
\newcommand{\MRhref}[1]{%
  \href{http://www.ams.org/mathscinet-getitem?mr=#1}{MR~#1}
}
\newcommand{\arxiv}[1]{%
  \href{http://www.arxiv.org/abs/#1}{arXiv:~#1}
}
%%%

%
%%%%%%%%%%%%%%%%%%%%%%%%%%%%%%%%%%%%%%%%%
%
% \input{glossary.tex}  %% Should be made by hand, using the file .glo. 
%
% \printindex           
%% This command does not work by itself. First, invoke the command
%% makeindex filename.idx (or simply: makeindex filename)
%% to produce filename.ind containing the index entries.
%% A Subsequent latex processing will print the index output in the
%% location of the \printindex command.
% 
% \bigskip
% % {\it
% {\footnotesize
% Author's Address: \medskip\\
% School of Mathematics \\
% Institute for Studies in Theoretical Physics and Mathematics (IPM) \\
% P.O. Box 19395--5746, Niavaran, Tehran, Iran \\
% E-mail: \href{mailto: safari@ipm.ir}{{\tt safari@ipm.ir}} \\
% Web: \href{http://math.ipm.ac.ir/safari/}{{\tt http://math.ipm.ac.ir/safari/}}
% }
% % }
%
%% Address
%
\bigskip
\noindent
\hypertarget{myaddress}%
{%
{\sc%
{\small%
\href{http://math.ipm.ac.ir/}{School of Mathematics}, 
\href{http://www.ipm.ac.ir/}{Institute for Studies in Theoretical Physics and Mathematics} (IPM),
P.O. Box 19395--5746, Niavaran, Tehran, Iran. \\
\href{mailto: safari@ipm.ir}{{\tt safari@ipm.ir}} \hfill
\href{http://math.ipm.ac.ir/safari/}{{\tt http://math.ipm.ac.ir/safari/}}
}
}
}
\end{document}